\documentclass[12pt]{amsart}
\usepackage{amsmath}
\usepackage{amsxtra}
\usepackage{amscd}
\usepackage{amsthm}
\usepackage{amsfonts}
\usepackage{amssymb}
\usepackage[mathscr]{eucal}
\usepackage{epsfig}
\usepackage{graphics}
\usepackage{accents}
\usepackage{physics}
\usepackage{hyperref}

\usepackage{mathtools}
\usepackage{tikz}
\usetikzlibrary{tikzmark,calc}
\usetikzlibrary{arrows,decorations.markings}
\usepackage[all]{xy}
\usepackage{here} 
\numberwithin{equation}{section}
\newtheorem{thm}{Theorem}[section]
\newtheorem{prop}[thm]{Proposition}
\newtheorem{lem}[thm]{Lemma}

\newtheorem{conj}[thm]{Conjecture}
\newtheorem{dfn}[thm]{Definition}

\newcommand{\C}{{\mathbb C}}
\newcommand{\Z}{{\mathbb Z}}

\newcommand{\Zz}{{\mathbb Z}^\times}

\newcommand{\E}{{\mathcal E}}

\newcommand{\gl}{\mathfrak{gl}}

\newcommand{\sln}{\mathfrak{sl}}

\newcommand{\La}{\Lambda}

\newcommand{\Sym}{\mathrm{Sym}}

\DeclareMathOperator\coker{coker}

\DeclareMathOperator\eva{ev}
\newcommand{\con}[1]{\left< #1\right>}
\newcommand{\Ug}{U_q\,\widehat{\mathfrak{gl}}_{m|n}}
\newcommand{\Uge}{\widetilde{U}_q\,\widehat{\mathfrak{gl}}_{m|n}}
\newcommand{\Us}{U_q\,\widehat{\mathfrak{sl}}_{m|n}}

\begin{document}

\begin{title}{Quantum toroidal algebra associated with $\gl_{m|n}$}
\end{title}
\author{Luan Bezerra and Evgeny Mukhin}

\address{LB: Department of Mathematics,
	Indiana University -- Purdue University -- Indianapolis,
	402 N. Blackford St., LD 270,
	Indianapolis, IN 46202, USA}\email{luanpere@iupui.edu}

\address{EM: Department of Mathematics,
	Indiana University -- Purdue University -- Indianapolis,
	402 N. Blackford St., LD 270,
	Indianapolis, IN 46202, USA}\email{emukhin@iupui.edu}

\begin{abstract}
	We introduce and study the quantum toroidal algebra $\E_{m|n}(q_1,q_2,q_3)$ associated with the  superalgebra $\gl_{m|n}$ with $m\neq n$, where the parameters satisfy $q_1q_2q_3=1$. 
	
	We give an evaluation map. The evaluation map is a surjective homomorphism of algebras $\E_{m|n}(q_1,q_2,q_3) \to \Uge$
	to the quantum affine algebra associated with the superalgebra $\gl_{m|n}$ at level $c$ completed with respect to the homogeneous grading, where $q_2=q^2$ and $q_3^{m-n}=c^2$.
	
	We also give a bosonic realization of level one $\E_{m|n}(q_1,q_2,q_3)$-modules.

\end{abstract}

\maketitle

\section{Introduction}
Quantum toroidal algebras were introduced in \cite{GKV} motivated by the study of Hecke operators in algebraic surfaces. 
Since that time, the quantum toroidal algebras, especially those associated with $\gl_n$, were found to have many applications in geometry, algebra, and mathematical physics. 

We list a few facts involving quantum toroidal algebras of type A.
The quantum toroidal algebras appear as Hall algebras of elliptic curves, \cite{BS}, \cite{SV1}, they also act on equivariant K-groups of Hilbert schemes and Laumon moduli spaces, \cite{FT},  \cite{SV2}, \cite{T1}. The quantum toroidal algebras are natural dual objects to double affine Hecke algebras, \cite{VV}. The quantum toroidal algebras provide integrable systems of XXY-type, among them is a deformation of quantum KdV flows, \cite{FJM1}. Characters of representations of quantum toroidal algebras appear in topological field theory, \cite{FJMM1}, AGT conjecture, \cite{AFS}. The full list is much longer.

\medskip

In this paper, we introduce the quantum toroidal algebras $\E_{m|n}(q_1,q_2,q_3)$ related to the superalgebras $\gl_{m|n}$, with $m\neq n$ and standard parity, and initiate their study.
In our mind, this subject is long overdue. We expect these algebras to have many properties similar to the quantum toroidal algebras $\E_{m|0}(q_1,q_2,q_3)$ associated with $\gl_m$ which can be used in similar way, but with various new features occurring due to the supersymmetry. In particular, our future goal is to study the corresponding integrable systems. 

\medskip

We start by introducing the algebras $\E_{m|n}(q_1,q_2,q_3)$ with $m\neq n$ and standard parity.  As in the even case, they depend on the complex parameters $q_1,q_2, q_3$ such that $q_1q_2q_3=1$. We require that the algebra $\E_{m|n}(q_1,q_2,q_3)$ has a ``vertical" quantum affine subalgebra $\Ug$ in the new Drinfeld realization, a ``horizontal'' quantum affine subalgebra $\Us$ given in Chevalley generators, and a symmetry with respect to the parity change, see \eqref{cph}. We always have $q_2=q^2$. This leads us to the generators and relations presentation of $\E_{m|n}(q_1,q_2,q_3)$, see Definition \ref{defE}. Naturally, the algebra $\E_{m|n}(q_1,q_2,q_3)$ is generated by currents $E_i(z), F_i(z)$, and half currents
$K_i^\pm(z)$, $i=0,\dots,m+n-1$, labeled by nodes of the affine Dynkin diagram of type $\widehat \sln_{m|n}$ and standard parity, see Figure \ref{fig}, and the relations are written in terms of the corresponding Cartan matrix. Similar to the even case, the quantum toroidal algebra $\E_{m|n}(q_1,q_2,q_3)$ has a two-dimensional center.

Then we note a few properties, including some isomorphisms and a topological Hopf algebra structure, see Sections \ref{secIso}, \ref{secH}. After that, we use bosonization techniques to construct level one representations of $\E_{m|n}(q_1,q_2,q_3)$, see Theorem \ref{T1}. Our formulas are built on work \cite{KSU} and generalize the known result in the even case \cite{S}.
We expect that the irreducible level one modules stay irreducible when restricted to the vertical $\Ug$ algebra. However, unlike the even case, the precise structure of irreducible level one modules for the quantum affine $\gl_{m|n}$ is not fully understood, see \cite{KSU}, \cite{K2}, and Conjecture \ref{conjM}. 

Finally, we proceed to the evaluation map. The evaluation map is a surjective algebra homomorphism $\E_{m|n}(q_1,q_2,q_3) \to \Uge$ to the quantum affine algebra at level $c$ completed with respect to the homogeneous grading, where $q_3^{m-n}=c^2$, see Theorem \ref{T2}. The evaluation map has the property that its restriction to the vertical subalgebra is the identity map. In the even case, the evaluation map was found in \cite{M2}, see also \cite{FJM2}. 

\medskip

Many results for $\E_{m|n}(q_1,q_2,q_3)$ need to be established, and we plan to address it in the follow-up papers. In particular, similar to the even case, we expect to obtain the Miki automorphism, see \cite{M1}, the shuffle algebra realization, see \cite{N}, the PBW type theorem, see \cite{T2}, the category $\mathcal O$, see \cite{FJMM1}, the fusion subalgebras, see \cite{FJMM2}, the integrable systems and Bethe ansatz, see \cite{FJMM3} with proper modifications.

In the supersymmetric case the Dynkin diagram of an algebra is not unique. It is straightforward to generalize our definition to other Dynkin diagrams associated with $\widehat\sln_{m|n}$. We expect that it gives the same algebra. For $m\neq n$ this is indeed so, see \cite{BM}. Moreover, one can generalize the definition of quantum toroidal algebras to other superalgebras in a similar fashion. In fact, quantum toroidal algebra related to the superalgebra D(2,1,$\alpha$) appeared implicitly in \cite{FJM3}. 

\medskip

The paper is organized as follows. 
In Section \ref{secdef} we define the quantum toroidal algebra associated with $\gl_{m|n}$ and give a few properties.
In Section \ref{secbos} we describe level 1 modules. In Section \ref{ev} we give the evaluation map. In Appendix \ref{app1} we collect some useful commutation relations. In Appendix \ref{Asl} we describe our conventions for the algebra $U_q\widehat{\sln}_{m|n}$.

\bigskip

{\bf Acknowledgments.} This work was partially supported by a grant from the Simons Foundation \#353831. L.B. was supported by the CNPq-Brazil grant 210375/2014-0. 
 
\section{Quantum toroidal \texorpdfstring{$\gl_{m|n}$}{glm|n}}{\label{secdef}}
Assume $m,n\geq 1$ and $m\neq n$. 

In this section we introduce the quantum toroidal algebra associated with $\gl_{m|n}$, denoted by $\E_{m|n}$, and collect a few properties.

\subsection{Definition of \texorpdfstring{$\E_{m|n}$}{Em|n}. }

Fix $d,q \in \mathbb{C}^\times$ and define 
$$q_1=d\,q^{-1},\quad q_2=q^2,\quad q_3=d^{-1}q^{-1}.$$

Note that $q_1q_2q_3=1$. Assume $q_1^{n_1}q_2^{n_2}q_3^{n_3}=1$, $n_1, n_2,n_3\in \mathbb{Z}$, iff  $n_1=n_2=n_3$. Fix also $d^{1/2}, q^{1/2}\in \mathbb{C}^\times$ such that $(d^{1/2})^2=d,\; (q^{1/2})^2=q$.

Let $I=\{1,\dots,m+n-1\}$. Let $\hat{I}^+=\{1,\dots,m-1\}$, $\hat{I}^-=\{m+1,\dots,m+n-1\}$, $\hat{I}^1=\{0,m\}$ and $\hat{I}=\hat{I}^+\cup \hat{I}^-\cup \hat{I}^1$. In particular, if $n=1$, we have $\hat{I}^-=\emptyset$, and if $m=1$, $\hat{I}^+=\emptyset$. The elements in these sets are to be read modulo $m+n$.

Let $\hat{A}=(A_{i,j})_{i,j\in \hat{I}}$ be the  Cartan matrix of $\widehat{\mathfrak{sl}}_{m|n}$ and $A=(A_{i,j})_{i,j\in I}$ be the  Cartan matrix of $\mathfrak{sl}_{m|n}$, both with the standard choice of parity. Namely, the odd simple roots correspond to $i\in \hat{I}^1=\{0,m\}$. We set $|i|=1$, $i \in \hat{I}^1$, and $|i|=0$ otherwise. We have $\det(A)=m-n\neq 0,$
${\rm {det}} (\hat A)= 0.$

For $i\in \hat{I}^+\cup \{m\}$, let  $M_{i-1,i}=-1, M_{i,i-1}=1$. For $j\in \hat{I}^-\cup \{0\}$, let 
$M_{j-1,j}=1, M_{j,j-1}=-1$. Let also  $M_{i,j}=0,$ if $i\neq j\pm 1$.
Define the matrix $\hat{M}=(M_{i,j})_{i,j\in \hat{I}}$. We have

\[\hat{A}=\begin{pmatrix}
	0 & -1 &&&  &&&1 \\
	-1 & 2 & -1   &&&&& \\
	&   -1& \ddots & \ddots  & \\ 
	& & \ddots & 2 & -1 \\
	& & &  -1 & 0 \tikzmark{k}& 1 &  \\
	& & &  & 1 & -2 &\ddots \\
	& & & &  & \ddots &\ddots & 1  \\ 
	1 &  &  &&&  & 1 & -2
	\end{pmatrix}
	,\quad  \hat{M}=\begin{pmatrix}
	0 & -1 &&&\tikzmark{k3}&&&-1 \\
	1 & 0 & -1   &&&&& \\
	&   1& \ddots & \ddots  & \\ 
	& & \ddots & 0 & -1 \\
	& & &  1 & 0\tikzmark{k2}& 1 &  \\
	& & &  & -1 & 0 &\ddots \\
	& & & &  & \ddots &\ddots & 1  \\ 
	1 &  &  &&&  & -1 & 0
	\end{pmatrix}.\]

	\begin{tikzpicture}[overlay,remember picture]
		\draw[<-,dashed] ($(pic cs:k)+(-0.1,-2)$) to
		($(pic cs:k)+(-0.1,-2.5)$) node[below] {$m$};
		\draw[<-,dashed] ($(pic cs:k2)+(-0.1,-2)$) to
		($(pic cs:k2)+(-0.1,-2.5)$) node[below] {$m$};
		\draw[<-,dashed] ($(pic cs:k2)+(3.1,0.10)$) to
		($(pic cs:k2)+(3.6,0.10)$) node[right] {$m$};
		\draw[<-,dashed] ($(pic cs:k3)+(3.2,0.10)$) to
		($(pic cs:k3)+(3.7,0.10)$) node[right] {$0$};
	\end{tikzpicture}
	
	\vspace{1cm}

	Note that $\hat{A}$ is symmetric and $ \hat{M}$ is skew-symmetric.
	
	\begin{dfn}{\label{defE}}
		The \textit{quantum toroidal algebra associated with $\gl_{m|n}$} is the unital associative superalgebra $\E_{m|n}=\E_{m|n}(q_1,q_2,q_3)$  generated by $E_{i,k},F_{i,k},H_{i,r}$, 
		and invertible elements $K_i$, $C$, where $i\in \hat{I}$, $k\in \Z$, $r\in\Zz$, subject to the relations \eqref{relCK}-\eqref{Serre8} below.
		The parity of the generators is given by $|E_{i,k}|=|F_{i,k}|=|i|$ and $0$ in all other cases.
	\end{dfn}

	The defining relations are given in terms of generating series
	\begin{align*}
		E_i(z)=\sum_{k\in\Z}E_{i,k}z^{-k}, \quad
		F_i(z)=\sum_{k\in\Z}F_{i,k}z^{-k}, \quad
		K^{\pm}_i(z)=K_i^{\pm1}\exp\Bigl(\pm(q-q^{-1})\sum_{r>0}H_{i,\pm r}z^{\mp r}\Bigr)\,. 
	\end{align*}
	
	We use the notation $[X,Y]_a=XY-(-1)^{|X||Y|}aYX.$ For simplicity, we write $[X,Y]_1=[X,Y]$.
	
	\medskip
	Let also $\displaystyle{\delta\left(z\right)=\sum_{n\in \mathbb{Z}} z^n}$.
	\medskip
	
	The relations are as follows.
	\bigskip

	\noindent{\bf $C,K$ relations}
	\begin{align}{\label{relCK}}
		  &\text{$C$ is central}, 
		  &&K_iK_j=K_jK_i,
		  &&K_iE_j(z)K_i^{-1}=q^{A_{i,j}}E_j(z),
		  &&K_iF_j(z)K_i^{-1}=q^{-A_{i,j}}F_j(z).
	\end{align}
	\bigskip
	
	\noindent{\bf $K$-$K$, $K$-$E$ and $K$-$F$ relations}
	\begin{align}
		  & K^\pm_i(z)K^\pm_j (w) = K^\pm_j(w)K^\pm_i (z)\,,                                                                          
		\label{KK1}\\
		  & \frac{d^{M_{i,j}}C^{-1}z-q^{A_{i,j}}w}{d^{M_{i,j}}Cz-q^{A_{i,j}}w}                                                        
		K^-_i(z)K^+_j (w) 
		=
		\frac{d^{M_{i,j}}q^{A_{i,j}}C^{-1}z-w}{d^{M_{i,j}}q^{A_{i,j}}Cz-w}
		K^+_j(w)K^-_i (z)\,,
		\label{KK2}\\
		  & (d^{M_{i,j}}z-q^{A_{i,j}}w)K_i^\pm(C^{-(1\pm1)/2}z)E_j(w)=(d^{M_{i,j}}q^{A_{i,j}}z-w)E_j(w)K_i^\pm(C^{-(1\pm1) /2}z)\,,   
		\label{KE}\\
		  & (d^{M_{i,j}}z-q^{-A_{i,j}}w)K_i^\pm(C^{-(1\mp1)/2}z)F_j(w)=(d^{M_{i,j}}q^{-A_{i,j}}z-w)F_j(w)K_i^\pm(C^{-(1\mp1) /2}z)\,. 
		\label{KF}
	\end{align}
	\bigskip

	\noindent{\bf $E$-$F$ relations}
	\begin{align}\label{EF}
		  & [E_i(z),F_j(w)]=\frac{\delta_{i,j}}{q-q^{-1}} 
		(\delta\left(C\frac{w}{z}\right)K_i^+(w)
		-\delta\left(C\frac{z}{w}\right)K_i^-(z))\,.
	\end{align}
	\bigskip
	
	\noindent{\bf $E$-$E$ and $F$-$F$ relations}
	\begin{align}
		  & [E_i(z),E_j(w)]=0\,, \quad [F_i(z),F_j(w)]=0\, \quad                                                   & (A_{i,j}=0)\,,    
		\\
		  & (d^{M_{i,j}}z-q^{A_{i,j}}w)E_i(z)E_j(w)=(-1)^{|i||j|}(d^{M_{i,j}}q^{A_{i,j}}z-w)E_j(w)E_i(z)\, \quad   & (A_{i,j}\neq0)\,, 
		\\
		  & (d^{M_{i,j}}z-q^{-A_{i,j}}w)F_i(z)F_j(w)=(-1)^{|i||j|}(d^{M_{i,j}}q^{-A_{i,j}}z-w)F_j(w)F_i(z)\, \quad & (A_{i,j}\neq0)\,. 
	\end{align}
	\bigskip
	
	\noindent{\bf Serre relations}
	\begin{align}
		  & \Sym_{{z_1,z_2}}[E_i(z_1),[E_i(z_2),E_{i\pm1}(w)]_q]_{q^{-1}}=0\,\quad & (i\not\in \hat{I}^1)\,,\label{Serre1} \\ 
		  & \Sym_{{z_1,z_2}}[F_i(z_1),[F_i(z_2),F_{i\pm1}(w)]_q]_{q^{-1}}=0\,\quad & (i\not\in \hat{I}^1)\,.               
		\label{Serre2}
	\end{align}
	If $mn\neq 2$,
	\begin{align}
		  & \Sym_{{z_1,z_2}}[E_i(z_1),[E_{i+1}(w_1),[E_i(z_2),E_{i-1}(w_2)]_q]_{q^{-1}}]=0\,\quad & (i\in \hat{I}^1)\,,\label{Serre3} 
		\\
		  & \Sym_{{z_1,z_2}}[F_i(z_1),[F_{i+1}(w_1),[F_i(z_2),F_{i-1}(w_2)]_q]_{q^{-1}}]=0\,\quad & (i\in \hat{I}^1)\,.\label{Serre4} 
	\end{align}
	
	If $(m,n)=(2,1)$,
	\begin{align}
		\label{Serre5} & \Sym_{{z_1,z_2}}\Sym_{{w_1,w_2}}[E_0(z_1),[E_2(w_1),[E_{0}(z_2),[E_2(w_2),E_{1}(y)]_q]]]_{q^{-1}}=    \\ \notag
		               & =\Sym_{{z_1,z_2}}\Sym_{{w_1,w_2}}[E_2(w_1),[E_0(z_1),[E_{2}(w_2),[E_0(z_2),E_{1}(y)]_q]]]_{q^{-1}}\,, \\
		\label{Serre6} & \Sym_{{z_1,z_2}}\Sym_{{w_1,w_2}}[F_0(z_1),[F_2(w_1),[F_{0}(z_2),[F_2(w_2),F_{1}(y)]_q]]]_{q^{-1}}=    \\ \notag
		               & =\Sym_{{z_1,z_2}}\Sym_{{w_1,w_2}}[F_2(w_1),[F_0(z_1),[F_{2}(w_2),[F_0(z_2),F_{1}(y)]_q]]]_{q^{-1}}\,. 
	\end{align}
	
	If $(m,n)=(1,2)$,
	\begin{align}
		\label{Serre7} & \Sym_{{z_1,z_2}}\Sym_{{w_1,w_2}}[E_0(z_1),[E_1(w_1),[E_{0}(z_2),[E_1(w_2),E_{2}(y)]_q]]]_{q^{-1}}=    \\ \notag
		               & =\Sym_{{z_1,z_2}}\Sym_{{w_1,w_2}}[E_1(w_1),[E_0(z_1),[E_{1}(w_2),[E_0(z_2),E_{2}(y)]_q]]]_{q^{-1}}\,, \\
		\label{Serre8} & \Sym_{{z_1,z_2}}\Sym_{{w_1,w_2}}[F_0(z_1),[F_1(w_1),[F_{0}(z_2),[F_1(w_2),F_{2}(y)]_q]]]_{q^{-1}}=    \\ \notag
		               & =\Sym_{{z_1,z_2}}\Sym_{{w_1,w_2}}[F_1(w_1),[F_0(z_1),[F_{1}(w_2),[F_0(z_2),F_{2}(y)]_q]]]_{q^{-1}}\,. 
	\end{align}

	Note that \begin{align*}
	K:=\prod_{i\in \hat{I}}K_i
	\end{align*} is a central element.
	
	The relations \eqref{KK1}-\eqref{KF} are equivalent to 
	
	\begin{align}
		  & [H_{i,r},E_j(z)]= \frac{[r A_{i,j}]}{r}d^{-rM_{i,j}}C^{-(r+|r|)/2}           
		\,z^r E_j(z)\,,
		\\
		  & [H_{i,r},F_j(z)]=- \frac{[r A_{i,j}]}{r}d^{-rM_{i,j}}C^{-(r-|r|)/2}          
		\,z^r F_j(z)\,,
		\\
		  & [H_{i,r},H_{j,s}]=\delta_{r+s,0} \cdot  \frac{[r A_{i,j}]}{r}d^{-rM_{i,j}}\, 
		\frac{C^r-C^{-r}}{q-q^{-1}}\,,\label{h_rel}
	\end{align}
	for all $r\in \Zz$, $i,j \in \hat{I}$, where $[k]=\frac{q^k-q^{-k}}{q-q^{-1}}$.
	
	\medskip
	
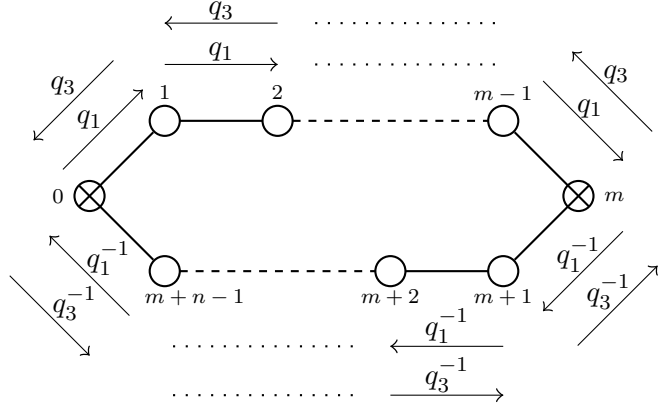
\begin{figure}[h]
	\centering
	\begin{tikzpicture}
		\draw[thick] (0,0) circle (0.2cm);
			\draw[rotate=45,thick] (-0.2,0)--(0.2,0);
			\draw[rotate=-45,thick] (-0.2,0)--(0.2,0);
		\draw[thick] (1,1) circle (0.2cm);
		\draw[thick] (2.5,1) circle (0.2cm);
		\draw[thick] (5.5,1) circle (0.2cm);
		\draw[thick] (6.5,0) circle (0.2cm);
		    \draw[shift={(6.5,0)}, rotate=45,thick] (-0.2,0)--(0.2,0);
			\draw[shift={(6.5,0)}, rotate=-45,thick] (-0.2,0)--(0.2,0);
		\draw[thick] (5.5,-1) circle (0.2cm);
		\draw[thick] (4,-1) circle (0.2cm);
		\draw[thick] (1,-1) circle (0.2cm);
		
		\draw[rotate=45,thick] (0.2,0)--(1.2,0);
		\draw[rotate=-45,thick] (0.2,0)--(1.2,0);
		\draw[shift={(6.5,0)}, rotate=135,thick] (0.2,0)--(1.2,0);
		\draw[shift={(6.5,0)}, rotate=-135,thick] (0.2,0)--(1.2,0);
		\draw[thick] (1.2,1)--(2.3,1);
		\draw[dashed,thick] (2.7,1)--(5.3,1);
		\draw[thick] (4.2,-1)--(5.3,-1);
		\draw[dashed,thick] (1.2,-1)--(3.8,-1);
		
		\node [below] at (1.4,-1.1) {\tiny $m+n-1$};
		\node [below] at (4,-1.1) {\tiny $m+2$};
		\node [below] at (5.5,-1.1) {\tiny $m+1$};
		\node [above] at (5.5,1.1) {\tiny $m-1$};
		\node [above] at (2.5,1.15) {\tiny $2$};
		\node [above] at (1,1.15) {\tiny $1$};
		\node [left] at (-0.2,0) {\tiny $0$};
		\node [right] at (6.7,0) {\tiny $m$};

		\draw [decoration={markings,mark=at position 1 with
    {\arrow[scale=1.5]{>}}},postaction={decorate}, shift={(1,1.75)}] (0,0) -- (1.5,0);
            \node [above] at (1.75,1.65) {\small $q_1$};
        \draw [decoration={markings,mark=at position 1 with
    {\arrow[scale=1.5]{>}}},postaction={decorate}, shift={(1,2.3)}]  (1.5,0)--(0,0);
            \node [above] at (1.75,2.2) {\small $q_3$};
            
        \draw [decoration={markings,mark=at position 1 with
    {\arrow[scale=1.5]{>}}},postaction={decorate} ,rotate=135, shift={(0.5,0)}] (0,0) -- (0,-1.5);
            \node [above] at (0, 0.75) {\small $q_1$};
        \draw [decoration={markings,mark=at position 1 with
    {\arrow[scale=1.5]{>}}},postaction={decorate}, rotate=135, shift={(1.05,0)}]  (0,-1.5)--(0,0);
            \node [above] at (-0.35,1.2) {\small $q_3$};
            
         \draw [decoration={markings,mark=at position 1 with
    {\arrow[scale=1.5]{>}}},postaction={decorate}, shift={(5.5,1)},rotate=45] (0.75,0) -- (0.75,-1.5);
            \node [above] at (6.65, 0.9) {\small $q_1$};
        \draw [decoration={markings,mark=at position 1 with
    {\arrow[scale=1.5]{>}}},postaction={decorate},  shift={(5.5,1)},rotate=45] (1.3,-1.5)--(1.3,0);
            \node [above] at (7,1.35) {\small $q_3$};

		\draw [decoration={markings,mark=at position 1 with
    {\arrow[scale=1.5]{>}}},postaction={decorate}, shift={(4,-2)}] (1.5,0)--(0,0);
            \node [above] at (4.75,-2.1) {\small $q_1^{-1}$};
        \draw [decoration={markings,mark=at position 1 with
    {\arrow[scale=1.5]{>}}},postaction={decorate}, shift={(4,-2.65)}]  (0,0)--(1.5,0);
            \node [above] at (4.75,-2.75) {\small $q_3^{-1}$};    
            
        \draw [decoration={markings,mark=at position 1 with
    {\arrow[scale=1.5]{>}}},postaction={decorate}, shift={(5.5,-1)},rotate=-45] (0.75,1.5)--(0.75,0);
            \node [above] at (6.5, -1.15) {\small $q_1^{-1}$};
        \draw [decoration={markings,mark=at position 1 with
    {\arrow[scale=1.5]{>}}},postaction={decorate},  shift={(5.5,-1)},rotate=-45] (1.4,0)--(1.4,1.5);
            \node [above] at (6.9,-1.65) {\small $q_3^{-1}$};    
            
        \draw [decoration={markings,mark=at position 1 with
    {\arrow[scale=1.5]{>}}},postaction={decorate}, shift={(0,0)},rotate=225] (0.75,1.5)--(0.75,0);
            \node [above] at (0.25, -1.2) {\small $q_1^{-1}$};
        \draw [decoration={markings,mark=at position 1 with
    {\arrow[scale=1.5]{>}}},postaction={decorate},  shift={(0,0)},rotate=225] (1.5,0)--(1.5,1.5);
            \node [above] at (-0.2,-1.8) {\small $q_3^{-1}$};        
            
        \draw[loosely dotted,thick] (3,1.75)--(5.5,1.75);
        \draw[loosely dotted,thick] (3,2.3)--(5.5,2.3);    
            
		\draw[loosely dotted,thick] (3.5,-2)--(1,-2);
        \draw[loosely dotted,thick] (3.5,-2.65)--(1,-2.65);  
	\end{tikzpicture}
	\caption{Standard Dynkin diagram of type $\widehat \sln_{m|n}$.}
	\label{fig}
\end{figure}

	The poles of the correlation functions of currents $E_i(z)$ are depicted in Figure \ref{fig}. For example, the correlation function of $E_0(z)E_1(w)$ has a pole at $z=q_1w$, while the correlation function of $E_1(z)E_0(w)$ has a pole at $z=q_3w$. The poles of the correlation functions of the currents $F_i(z)$ are obtained from the Figure \ref{fig} replacing $q$ by $q^{-1}$, i.e., $q_1^{\pm 1}$ is replaced by $q_3^{\mp 1}$, and $q_3^{\pm 1}$ by $q_1^{\mp 1}$.
	
	\bigskip
	
	\subsection{Horizontal and vertical subalgebras}
	
	In this section, we define the horizontal and vertical subalgebras of $\E_{m|n}$.
	
	See Appendix \ref{Asl} for the definitions of the quantum affine algebras $U_q\,\widehat{\mathfrak{sl}}_{m|n}$ and $U_q\,\widehat{\mathfrak{gl}}_{m|n}$.
	
	We denote the subalgebra of $\E_{m|n}$ generated by $E_i(z), F_i(z), K_i^\pm (z), C, i \in I$,  by $U^{ver}_q\widehat{\mathfrak{sl}}_{m|n}$ and we call it the {\it vertical quantum affine $\mathfrak{sl}_{m|n}$}. If $K_0^\pm(z)$ are also included, the resulting subalgebra is denoted by $U^{ver}_q\widehat{\mathfrak{gl}}_{m|n}$ and called the {\it vertical quantum affine $\gl_{m|n}$}. Note that $U^{ver}_q\widehat{\mathfrak{sl}}_{m|n}$, $U^{ver}_q\widehat{\mathfrak{gl}}_{m|n}$ are given in new Drinfeld realization.
	
	The currents $K_0^\pm(z)$ do not commute with $U^{ver}_q\widehat{\mathfrak{sl}}_{m|n}$. To obtain a current in $U^{ver}_q\widehat{\mathfrak{gl}}_{m|n}$ commuting with $U^{ver}_q\widehat{\mathfrak{sl}}_{m|n}$ we proceed as follows. 
	
	For each $r\in \Zz$, $\det([r A_{i,j}]d^{-rM_{i,j}})_{i,j \in \hat{I}}=[r]^{m+n}\left(d^{r(m-n)}+d^{r(n-m)}-q^{r(m-n)}-q^{r(n-m)}\right)\neq 0$. Thus, the system
	\begin{align}{\label{sys1}}
		\sum_{i\in \hat{I}}\gamma_{i,r}[rA_{i,j}]d^{-rM_{i,j}}=0 \quad (j\in I)\,, 
	\end{align}
	has a one-dimensional space of solutions. The element $H^{ver}_r= \sum_{i\in \hat{I}}\gamma_{i,r}H_{i,r} \in U^{ver}_q\widehat{\mathfrak{gl}}_{m|n}$ commutes with $U^{ver}_q\widehat{\mathfrak{sl}}_{m|n}\subset U^{ver}_q\widehat{\mathfrak{gl}}_{m|n}$. Such element is unique up to scalar. We fix a normalization by requiring $\gamma_{0,r}=1,\, r\in \Z_{<0}$, and
	\begin{align*}
		  [H^{ver}_r,H^{ver}_{s}]=\delta_{r+s,0}[(n-m)r]\frac{1}{r}\frac{C^r-C^{-r}}{q-q^{-1}}. 
	\end{align*}
	
	Set $H^{ver}(z)=\sum_{r\in \Z^\times} H^{ver}_r z^{-r}$.

	\begin{lem}{\label{lemv}}
		The subalgebra $U^{ver}_q\widehat{\mathfrak{gl}}_{m|n}$ is isomorphic to $\Ug$ for generic values of parameters.
		
		\begin{proof}
			We have a homomorphism $v:\Ug \rightarrow \E_{m|n}$ given by
			\begin{align*}
				  & x^+_i(z)\mapsto E_i(d^{-i}z), \quad  x^-_i(z)\mapsto F_i(d^{-i}z),\quad k_i^\pm(z)\mapsto K_i^\pm(d^{-i}z)\quad                      & (i\in \hat{I}^+\cup \{m\})\,, \\ 
				  & x^+_{j}(z)\mapsto E_{j}(d^{-2m+j}z), \quad  x^-_{j}(z)\mapsto F_{j}(d^{-2m+j}z),\quad k_{j}^\pm(z)\mapsto K_{j}^\pm(d^{-2m+j}z)\quad & (j\in \hat{I}^-)\,,           \\ 
				&c\mapsto C,\quad h(z)\mapsto H^{ver}(z)\,.
			\end{align*}
			
			The evaluation map constructed in Theorem \ref{T2} produces a left-inverse of $v$, for generic values of parameters, see Lemma \ref{lem1}. Thus, $v$ is an embedding with image $U^{ver}_q\widehat{\mathfrak{gl}}_{m|n}$.
		\end{proof}
	\end{lem}
	
	We denote the subalgebra of $\E_{m|n}$ generated by $E_{i,0}, F_{i,0}, K_i, i \in \hat{I}$, by $U^{hor}_q\widehat{\mathfrak{sl}}_{m|n}$ and we call it the {\it horizontal quantum affine $\mathfrak{sl}_{m|n}$}. Note that $U^{hor}_q\widehat{\mathfrak{sl}}_{m|n}$ is given in Drinfeld-Jimbo realization.
	
	We have a homomorphism $h:\Us \rightarrow \E_{m|n}$ given by
	\begin{align*}
		  & e_i\mapsto E_{i,0}\,, \quad  f_i\mapsto F_{i,0}\,,\quad t_i\mapsto K_i & (i\in \hat{I})\,, 
	\end{align*}    
	and its image is $U^{hor}_q\widehat{\mathfrak{sl}}_{m|n}$. 
	
	\begin{conj}
		The homomorphism $h:\Us \rightarrow \E_{m|n}$  
		is injective. In particular, $U^{hor}_q\widehat{\mathfrak{sl}}_{m|n}$ is isomorphic to $\Us$.
	\end{conj}
	This conjecture is proved in \cite{BM} for $m\neq n$, $m+n>3$, and generic values of parameters. 
	
	\bigskip
	\subsection{Isomorphisms}{\label{secIso}} In this section, we list some isomorphisms involving superalgebra $\E_{m|n}$. In all cases, it is easy to check that the maps are even, invertible and respect the defining relations.

	For $a\in \mathbb{C}^\times$, the {\it shift of spectral parameter}  
	\begin{align*}
		s_a:\E_{m|n}\rightarrow \E_{m|n}, 
	\end{align*}
	is defined by
	\begin{align*}
		  & s_a(E_i(z))=E_i(az),\quad s_a(F_i(z))=F_i(az),\quad s_a(K_i^\pm(z))=K_i^\pm(az), \quad s_a(C)=C \quad & (i\in \hat{I}). 
	\end{align*}

	For each $j\in \hat{I}$, we have
	\begin{align*}
		\chi_j:\E_{m|n}\rightarrow \E_{m|n}, 
	\end{align*}
	defined by
	\begin{align*}
		  \chi_j(E_i(z))=E_i(z)z^{-\delta_{i,j}},\quad \chi_j(F_i(z))=F_i(z)z^{\delta_{i,j}},\quad \chi_j(K_i^\pm(z))=C^{\mp \delta_{i,j}}K_i^\pm(z),\quad \chi_j(C)=C\quad  (i\in \hat{I}). 
	\end{align*}
	
	\medskip
	The following isomorphisms change the parameters of the algebra.
	\medskip
	
	The {\it diagram isomorphism}
	\begin{align}{\label{diagiso}}
		  & \sigma:\E_{m|n}(q_1,q_2,q_3) \rightarrow \E_{m|n}(q_3,q_2,q_1), 
	\end{align}
	defined by
	\begin{align*}
		  & \sigma(E_i(z))=E_{m-i}(z),\quad  \sigma(F_i(z))=F_{m-i}(z),\quad  \sigma(K_i^\pm(z))=K_{m-i}^\pm(z), \quad \sigma(C)=C\quad & (i\in \hat{I}), 
	\end{align*}
	changes $d$ to $d^{-1}$.
	
	The {\it change of parity} isomorphism
	\begin{align}{\label{cph}}
		  & \tau:\E_{m|n}(q_1,q_2,q_3) \rightarrow \E_{n|m}(q_3^{-1},q_2^{-1},q_1^{-1}), 
	\end{align}
	defined by
	\begin{align*}
		  & \tau(E_i(z))=E_{-i}(z),\quad  \tau(F_i(z))=F_{-i}(z),\quad  \tau(K_i^\pm(z))=-K_{-i}^\pm(z), \quad \tau(C)=C\quad & (i\in \hat{I}), 
	\end{align*}
	changes $q$ to $q^{-1}$.
	
	We have $\sigma^2=\tau^2=Id$.
	
	\bigskip
	\subsection{Hopf superalgebra structure}{\label{secH}}
	
	\begin{prop}
		The superalgebra $\E_{m|n}$ has a topological Hopf superalgebra structure given on generators by
	\begin{align*}
		  & \Delta E_i(z)=E_i(z)\otimes 1 + K^-_i(z)\otimes E_i(C_1z),                                 \\
		  & \Delta F_i(z)=F_i(C_2z)\otimes K^+_i(z) + 1\otimes F_i(z),                                 \\
		  & \Delta K^+_i(z)=K^+_i(C_2z)\otimes K^+_i(z),                                               \\
		  & \Delta K^-_i(z)=K^-_i(z)\otimes K^-_i(C_1z),                                               \\
		  & \Delta C=C\otimes C,                                                                       \\
		  & \varepsilon(E_i(z))=\varepsilon(F_i(z))=0,\quad \varepsilon(K^\pm _i(z))=\varepsilon(C)=1, \\
		  & S(E_i(z))=-\left(K^-_i(C^{-1}z)\right)^{-1}E_i(C^{-1}z),                                   \\
		  & S(F_i(z))=-F_i(C^{-1}z)\left(K^+_i(C^{-1}z)\right)^{-1},                                   \\
		  & S(K^\pm_i(z))=\left(K^\pm _i(C^{-1}z) \right)^{-1},\quad S(C)=C^{-1},                      
	\end{align*}
	where $C_1=C\otimes 1$, $C_2=1\otimes C$. The maps $\Delta$ and $\varepsilon$ are extended to algebra homomorphisms, and the map $S$ to a superalgebra anti-homomorphism, $S(xy)=(-1)^{|x||y|}S(y)S(x)$.
	\begin{proof} The proof is done by straightforward computations. In the $U_q\widehat{\mathfrak{gl}}_{m|n}$ case, these formulas appeared in \cite{Z}. For $U_q\widehat{\mathfrak{gl}}_{m}$, a proof is given in \cite{DI}.
\end{proof}
	
	\end{prop}
	
 Note that the tensor product multiplication is defined for homogeneous elements $x_1, x_2, y_1, y_2 \in\E_{m|n} $ by $(x_1\otimes y_1)(x_2\otimes y_2)=(-1)^{|y_1||x_2|}x_1x_2\otimes y_1y_2$ and extended to the whole algebra by linearity.
	
	The vertical subalgebras $U^{ver}_q\widehat{\mathfrak{sl}}_{m|n}$ and $U^{ver}_q\widehat{\mathfrak{gl}}_{m|n}$ are Hopf subalgebras of $\E_{m|n}$.
	
	\bigskip
	\subsection{Grading}{\label{secG}}
	
	For each $i \in \hat{I}$, the superalgebra $\E_{m|n}$ has a $\Z$-grading given by 
	\begin{align*}
		&\deg_i(E_{j,k})=\delta_{i,j},\quad \deg_i(F_{j,k})=-\delta_{i,j},\\ &\deg_i(H_{j,r})=\deg_i  K_j=\deg_i(C)=0\quad  (j\in \hat{I},  k\in \Z, r\in \Zz).
	\end{align*}

	There is also the \textit{homogeneous $\Z$-grading} given by
	\begin{align*}
		\deg_{\delta}(E_{j,k})=\deg_{\delta}(F_{j,k})=k, \quad \deg_{\delta}(H_{j,r})=r, \quad \deg_{\delta}(K_j)=\deg_{\delta}(C)=0 \quad  (j\in \hat{I},  k\in \Z, r\in \Zz). 
	\end{align*}
	
	Thus the superalgebra $\E_{m|n}$ has a $\Z^{m+n+1}$-grading given by
	\begin{align*}
	    \deg (X)=\left(\deg_0(X), \deg_1(X),\dots,\deg_{m+n-1}(X);\deg_{\delta}(X) \right), \qquad X\in \E_{m|n}.
	\end{align*}
	
	We call an $\E_{m|n}$-module $V$ \textit{admissible} if for any $v\in V$ there exists an integer $N=N_v$ such that $Xv=0$ for all $X \in \E_{m|n}$ with $\deg_{\delta}(X)>N$.

	We say an $\E_{m|n}$-module has level $(k^{ver},k^{hor})$ if it has level $k^{ver}$ as a  $U^{ver}_q\widehat{\mathfrak{sl}}_{m|n}$-module, and level $k^{hor}$ as a  $U^{hor}_q\widehat{\mathfrak{sl}}_{m|n}$-module, i.e., if $(C,K)$ acts as $(q^{k^{ver}},q^{k^{hor}})$.

	\bigskip
	\section{Level (1,0) modules, bosonic picture}{\label{secbos}}
	In this section, we construct $\E_{m|n}$-modules of level $(1,0)$ using vertex operators.
	
	\subsection{Heisenberg algebra}
	
	Let $\mathcal{H}$ be the associative algebra generated by $H_{i,r}$, $c_{j,r}$,  $i\in \hat{I}$, $j \in \hat{I}^-\cup \{m\}$, $r\in \Zz$, satisfying
	\begin{align}
		  & [H_{i,r},H_{j,s}]=\delta_{r+s,0} \cdot  \frac{[r A_{i,j}][r]}{r}d^{-rM_{i,j}}, {\label{h_rel1}} \\
		  & [c_{i,r},c_{j,s}]=\delta_{i,j}\delta_{r+s,0} \cdot \dfrac{[r]^2}{r},\notag                             \\
		  & [H_{i,r},c_{j,s}]=0.\notag                                                                      
	\end{align}
	Note that \eqref{h_rel1} is equivalent to equation \eqref{h_rel} with $C=q$.
	
	Denote by $\mathcal{H}^\pm$ the (commutative) subalgebra generated by $H_{i,r},c_{j,r}$ with $\pm r>0$, $i \in \hat{I}$, $j\in \hat{I}^-\cup \{m\}$.
	
	Let $\mathcal{F}$ be the Fock space generated by a vector $v_0$ satisfying $H_{i,r}v_0=c_{j,r}v_0=0$, for $ r>0$, $i \in \hat{I}$, $j\in \hat{I}^-\cup \{m\}$. Thus, $\mathcal{F}$ is a free $\mathcal{H}^-$-module of rank 1
	$$\mathcal{F}=\mathcal{H}v_0=\mathcal{H}^-v_0.$$
	Moreover, since $\det([r A_{i,j}]d^{-rM_{i,j}})_{i,j \in \hat{I}}\neq 0$, $\mathcal{F}$ is an irreducible $\mathcal{H}$-module.
	
	\subsection{Level (1,0)  \texorpdfstring{$\E_{m|n}$}{Em|n}-modules}
	
	Let $ Q_{m|n}$ be the $\mathfrak{sl}_{m|n}$  root lattice and let $\mathbb{C}\{Q_{m|n}\}$ be a twisted group algebra of ${Q}_{m|n}$ generated by invertible elements $e^{\bar{\alpha}_i}, i\in I,$ satisfying the relations
	$$e^{\bar{\alpha}_i}e^{\bar{\alpha}_j}=\begin{cases}(-1)^{\braket{\bar{\alpha}_i}{\bar{\alpha}_j}} e^{\bar{\alpha}_j}e^{\bar{\alpha}_i} &(i,j\in \hat{I}^+\cup\{m\}),\\ e^{\bar{\alpha}_j}e^{\bar{\alpha}_i} &(i\, \text{or}\, j \in \hat{I}^-).\end{cases}$$
	
	Define $\varepsilon:\hat{I}\times \hat{I}\rightarrow \{\pm 1\}$ by
	$$\varepsilon(i,j)=\begin{cases}(-1)^{\braket{\bar{\alpha}_i}{\bar{\alpha}_j}} &(i,j\in \hat{I}^+\cup\{m\},\, i>j),\\
	(-1)^{1+\delta_{m,1}} &(i=0,\, j\in \{1,m\}),\\
	1 &(\text{otherwise}).\end{cases}$$
Note that 	
	$$e^{\bar{\alpha}_i}e^{\bar{\alpha}_j}=\varepsilon(i,j)\varepsilon(j,i)e^{\bar{\alpha}_j}e^{\bar{\alpha}_i} \qquad (i,j \in \hat{I}).$$ 
	Let $Q_c$ be the integral lattice generated by elements $c_{i},\, i\in \hat{I}^-\cup \{m\},$ with bilinear form given by
	\begin{align*}
		  & \braket{c_{i}}{c_{j}}=\delta_{i,j}. 
	\end{align*}
	
	Define $Q=Q_{m|n}\oplus Q_c$ and extend the bilinear forms on $Q_{m|n}$ and $Q_c$ to $Q$ by requiring $\braket{\bar{\alpha}_i}{c_{j}}=0$. Set also $\braket{\bar{\Lambda}}{c_{j}}=0$, for any $\mathfrak{gl}_{m|n}$ weight $\bar{\Lambda}$.

	Let $\mathbb{C}[Q_c]$ be the (commutative) group algebra of $Q_c$ and define $\mathbb{C}\{Q\}=\mathbb{C}\{Q_{m|n}\}\otimes \mathbb{C}[Q_c]$.
	
	For $\alpha=\sum_{l\in I} r_{l}\bar{\alpha}_{l}+\sum_{k\in \hat{I}^-\cup \{m\}}s_kc_k\in Q $, define  
	\begin{align}{\label{alpha decomp}}
	    e^{\alpha}=(e^{\bar{\alpha}_{1}})^{r_{1}}\cdots (e^{\bar{\alpha}_{m+n-1}})^{r_{m+n-1}}(e^{c_{m}})^{s_{m}}\cdots (e^{c_{m+n-1}})^{s_{m+n-1}}.
	\end{align}
	Then, $\{e^{\alpha}, \alpha \in Q\}$ is a basis of $\mathbb{C}\{Q\}$.
	
	Let $\Tilde{Q}\subset Q$ be the sublattice of rank $m+n-1$ generated by $\bar{\alpha}_i, \bar{\alpha}_m + c_m, \bar{\alpha}_j + c_j - c_{j-1}$, $i\in \hat{I}^+, j \in \hat{I}^-$, and let $\mathbb{C}\{\Tilde{Q}\}$ be the subalgebra of $\mathbb{C}\{Q\}$ spanned by $e^{\alpha}$, $\alpha\in \Tilde{Q}$.
	
	Following \cite{KW}, a $\widehat{\mathfrak{sl}}_{m|n}$ weight $\Lambda$ is a level $1$ partially integrable weight if and only if $\Lambda=\Lambda_i,\, i\not \in \hat{I}^1$, or $\Lambda=(1-a)\Lambda_0+a\Lambda_m,\, a\in \mathbb{C}$.
	
	Set
	\begin{align*}
		  & \Tilde{\Lambda}=\bar{\Lambda}_i                         & (\Lambda=\Lambda_i,\,i\in \hat{I}^+),                  \\
		  & \Tilde{\Lambda}=\bar{\Lambda}_j-\sum_{i=j}^{m+n-1}c_i   & (\Lambda=\Lambda_j,\,j\in \hat{I}^-),                  \\
		  & \Tilde{\Lambda}=a\bar{\Lambda}_m-a\sum_{i=m}^{m+n-1}c_i & (\Lambda=(1-a)\Lambda_0+a\Lambda_m,\,a\in \mathbb{C}). 
	\end{align*}
	
	Given a level $1$ partially integrable weight $\Lambda$, define the vector superspace
	\begin{align*}
		  & \mathcal{F}_{\Lambda}:=\mathcal{F}\otimes\mathbb{C}\{\Tilde{Q}\}e^{\Tilde{\Lambda}}. 
	\end{align*}
For $v\in \mathcal{F},\, \alpha \in \Tilde{Q}$, the parity of $v\otimes e^\alpha e^{\Tilde{\Lambda}}\in \mathcal{F}_{\Lambda}$ is $|v\otimes e^\alpha e^{\Tilde{\Lambda}}|=(1-(-1)^{r_m})/2$, where $r_m$ is the multiplicity of $\bar{\alpha}_m$ in $\alpha$ as in \eqref{alpha decomp}.

	Define an action of the algebras $ \mathcal{H}$ and $\mathbb{C}\{\Tilde{Q}\}$ on $\mathcal{F}_{\Lambda}$ as follows.
	
	For $v\in \mathcal{F},\, \alpha \in \Tilde{Q}$, set
	\begin{align*}
		  & x(v\otimes e^\alpha e^{\Tilde{\Lambda}})=(xv)\otimes e^\alpha e^{\Tilde{\Lambda}}              & (x\in \mathcal{H}),    \\
		  & e^\beta(v\otimes e^\alpha e^{\Tilde{\Lambda}})=v\otimes (e^\beta e^\alpha e^{\Tilde{\Lambda}}) & (\beta \in \Tilde{Q}). 
	\end{align*}
	In particular, $\mathcal{F}_{\Lambda}$ is a free $\mathcal{H}^- \otimes \mathbb{C}\{\Tilde{Q}\}$-module of rank $1$.
	
	Introduce the zero-mode linear operators $z^{\pm H_{i,0}},q^{\pm \alpha_{i,0}},z^{\pm c_{j,0}},\;i \in \hat{I},\; j\in \hat{I}^-\cup \{m\} $, acting on $\mathcal{F}_{\Lambda}$ as follows.
	
	For $v\otimes  e^\alpha e^{\Tilde{\Lambda}} \in \mathcal{F}_{\Lambda}$, with $\alpha=\sum_{l\in I} r_{l}\bar{\alpha}_{l}+\sum_{k\in \hat{I}^-\cup \{m\}}s_kc_k$, set
	\begin{align*}
		  & z^{\pm H_{i,0}}(v\otimes e^{\alpha}e^{\tilde{\Lambda}})=z^{\pm \braket{\bar{\alpha}_i}{\alpha+{\tilde{\Lambda}}}}d^{\pm \frac{1}{2}\sum_{l\in I} r_lA_{i,l}M_{i,l}}v\otimes  e^{\alpha}e^{\tilde{\Lambda}}, \\
		  & q^{\pm \alpha_{i,0}}(v\otimes e^\alpha e^{\tilde{\Lambda}})=q^{\pm \braket{\bar{\alpha}_i}{\alpha+{\tilde{\Lambda}}}}v\otimes  e^\alpha e^{\tilde{\Lambda}},                                                \\
		  & z^{\pm c_{j,0}}(v\otimes e^\alpha e^{\tilde{\Lambda}})=z^{\pm \braket{c_{j}}{\alpha+{\tilde{\Lambda}}}}v\otimes  e^\alpha e^{\tilde{\Lambda}}.                                                              
	\end{align*}
	
	For $i \in \hat{I}$ and $j\in \hat{I}^-\cup \{m\}$, let
	\begin{align*}
		  & H_i^\pm(z)=\sum_{r>0}\dfrac{H_{i,\pm r}}{[r]}z^{\mp r},    \\
		  & c_{j}^\pm(z)=\sum_{r> 0}\dfrac{c_{j,\pm r}}{[r]}z^{\mp r}, 
	\end{align*}
	and define
	\begin{align*}
		  & \Gamma^+_i(z)=z\exp( H_i^-(q^{-1}z))\exp(-H_i^+(z)) e^{\bar{\alpha}_i}z^{H_{i,0}},   \\
		  & \Gamma^-_i(z)=z\exp( -H_i^-(z))\exp(H_i^+(q^{-1}z)) e^{-\bar{\alpha}_i}z^{-H_{i,0}}, \\
		  & C^\pm_{j}(z)=\exp(\pm c^-_j(z))\exp(\mp c^+_j(z))e^{\pm c_{j}}z^{\pm c_{j,0}}.       
	\end{align*}
	
	Note that these currents act on the larger space $\mathcal{F}\otimes\mathbb{C}\{Q\}e^{\Tilde{\Lambda}}$. However, their products considered on Theorem \ref{T1} below preserve the subspace $\mathcal{F}_{\Lambda}\subset \mathcal{F}\otimes\mathbb{C}\{Q\}e^{\Tilde{\Lambda}} $.

	The following is proved by a direct computation.
	
	\begin{lem}\label{lem_h}
		For $i,j \in \hat{I}$, $r\in \Zz$, we have
		\begin{align}
			  & [H_{i,r},\Gamma^\pm_j(z) ]=\pm\dfrac{[r A_{i,j}]}{r}d^{-rM_{i,j}}q^{-(r\pm |r|)/2}z^r\,\Gamma^\pm_j(z). 
		\end{align}
		\qed
	\end{lem}
	
	Define the normal ordering by
	\begin{align*}
		& :x_r y_s:=\begin{cases}
		x_r y_s\; &(r<0)\\
		y_s x_r\; &(r\geq 0)\\
		\end{cases} & (x_r,y_r\in \{ H_{i,r}, c_{i,r}\}),\\
		  & :Ae^\alpha :=:e^\alpha  A:=e^\alpha  A                                                    & ( \alpha  \in Q, A \in \{z^{\pm H_{i,0}},q^{\pm \alpha_{i,0}},z^{\pm c_{j,0}}\}), \\
		  & :e^{\bar{\alpha}_i} e^{\bar{\alpha}_j}:=\varepsilon(i,j)e^{\bar{\alpha}_i+\bar{\alpha}_j} & (i,j \in \hat{I}),                                                                
	\end{align*}
	and extended inductively from right to left on larger products, e.g., $:abc:=:a(:bc:):$. 
	
	Given two currents $X(z), Y(w)$, we say that $X(z)Y(w)$ has contraction $\con{X(z)Y(w)}$ if $$X(z)Y(w)=\con{X(z)Y(w)}:X(z)Y(w):.$$
	
	In this text, all contractions $\con{X(z)Y(w)}$ are Laurent series converging to rational functions in the region  $|z|\gg|w|$.
	\begin{lem}
		For $i,j\in \hat{I}$, $k,l\in \hat{I}^-\cup\{m\}$ we have
		\begin{align}
			&\con{\Gamma^\pm_i(z)\Gamma^\pm_i(w)}=\left((z-w)(z-q^{\mp 2}w) \right)^{A_{i,i}/2},\label{gii} \\
			  & \con{\Gamma^\pm_i(z)\Gamma^\pm_j(w)}=\varepsilon(i,j)\left(z-d^{-M_{i,j}}q^{\mp 1}w\right)^{A_{i,j}}d^{A_{i,j}M_{i,j}/2} & (i\neq j),\label{gij}   \\
			&\con{\Gamma^\pm_i(z)\Gamma^\mp_i(w)}=\left((z-qw)(z-q^{-1}w) \right)^{-A_{i,i}/2},\label{giim} \\
			  & \con{\Gamma^\pm_i(z)\Gamma^\mp_j(w)}=\varepsilon(i,j)\left(z-d^{-M_{i,j}}w\right)^{-A_{i,j}}d^{-A_{i,j}M_{i,j}/2}        & (i\neq j), \label{gijm} \\
			&\con{C^\pm_{k}(z)C^\pm_{l}(w)}=\left(z-w\right)^{ \delta_{k,l}}\label{cc1},\\
			&\con{C^\pm_{k}(z)C^\mp_{l}(w)}=\left(z-w\right)^{- \delta_{k,l}}\label{cc2}.
		\end{align}
		\begin{proof}
			Let $\epsilon \in \{\pm 1\}$.
			
			Equations \eqref{gii}-\eqref{gijm} follow from 
			\begin{align*}
				&\con{z^{\pm \epsilon H_{i,0}}e^{\epsilon \bar{\alpha}_j}}=z^{\pm A_{i,j}}d^{\pm A_{i,j}M_{i,j}/2},  \\
				&\con{\exp(\pm \epsilon H_i^+(z))\exp(\epsilon H_i^-(w))}=\left(1-q\dfrac{w}{z}\right)^{\mp A_{i,i}/2}\left(1-q^{-1}\dfrac{w}{z}\right)^{\mp A_{i,i}/2}, \\
				  & \con{\exp(\pm \epsilon H_i^+(z))\exp(\epsilon H_j^-(w))}= \left(1-d^{-M_{ij}}\dfrac{w}{z}\right)^{\mp A_{i,j}} & (i\neq j). 
			\end{align*}
			The equations \eqref{cc1} and \eqref{cc2} follow from     
			\begin{align*}
				  & \con{z^{\pm \epsilon c_{k,0}}e^{\epsilon c_l}}=z^{\pm \delta_{k,l}},                                    \\
				  & \con{\exp(\pm \epsilon c_k^+(z))\exp(\epsilon c_l^-(w))}=\left(1-\frac{w}{z}\right)^{\mp \delta_{k,l}}. 
			\end{align*}
			These contractions are checked by a straightforward computation.
		\end{proof}
	\end{lem}
	
	Let $\partial_z$ be the $q$-difference operator
	$$\partial_zf(z)=\dfrac{f(qz)-f(q^{-1}z)}{(q-q^{-1})z}.$$
	
	\begin{thm}{\label{T1}}
		The following expressions define a graded admissible $\E_{m|n}-$module structure of level $(1,0)$ on $\mathcal{F}_{\Lambda}$. 
		\begin{align*}
			  & C=q, \quad K_i^{\pm1}= q^{\pm \alpha_{i,0}}, \quad H_{i,r}=H_{i,r}                                 & (i\in \hat{I}),    \\
			  & E_i(z)=\Gamma^+_i(z)\;                                                                             & (i\in \hat{I}^+),  \\
			  & F_i(z)=\Gamma^-_i(z)\;                                                                             & (i\in \hat{I}^+),  \\
			&E_m(z)=d^m\Gamma^+_m(z)C^+_{m}(d^m z)\, ,\\
			&F_m(z)=\Gamma^-_m(z)\,\partial_{z}[C^-_{m}(d^m z)]\, ,\\
			  & E_i(z)=d^{(2m-i)}\Gamma^+_i(z) :C^+_{i}(d^{2m-i} z)\,\partial_{z}[C^-_{i-1}(d^{2m-i} z)]: & (i\in \hat{I}^-),  \\
			  & F_i(z)=d^{(2m-i)}\Gamma^-_i(z) :C^+_{i-1}(d^{2m-i} z)\,\partial_{z}[C^-_{i}(d^{2m-i} z)]: & (i\in \hat{I}^-) , \\
			&E_0(z)=\Gamma^+_0(z)\,\partial_z[C^-_{m+n-1}(d^{m-n} z)]\, ,\\
			&F_0(z)=d^{m-n} \Gamma^-_0(z)C^+_{m+n-1}(d^{m-n} z) .
		\end{align*}
		\begin{proof}
			The $C, K$ relations are clear.
			
			The $H$-$H$, $H$-$E$ and $H$-$F$ relations follow from Lemma \ref{lem_h}. Note that $H_{i,r}$ commutes with $c_j(z)$ for all possible $i,j$.
			
			We now check the $E$-$E$ relations.
			
			If $i=j \in \hat{I}^+$, it follows from \eqref{gii} that
			\begin{align*}
				\Gamma^+_i(z)\Gamma^+_i(w)=\Gamma^+_i(w)\Gamma^+_i(z) \left(\frac{z-q^{-2}w}{q^{-2}z-w}\right)^{A_{i,i}/2}, 
			\end{align*}
			which is equivalent to the $E$-$E$ relation.
			
			If $i\in \hat{I}^+$, $j\in \hat{I}$ and $i\neq j$ we use \eqref{gij}  to get  
			\begin{align*}
				\Gamma^+_i(z)\Gamma^+_j(w)=\Gamma^+_j(w)\Gamma^+_i(z) \left(\frac{d^{M_{i,j}}z-q^{-1}w}{w-d^{M_{i,j}}q^{-1}z}\right)^{A_{i,j}}\frac{\varepsilon(i,j)}{\varepsilon(j,i)}, 
			\end{align*}
			but in this case $\varepsilon(i,j)=(-1)^{A_{i,j}}\varepsilon(j,i)$, which is the needed sign.
			
			The cases with $i \in \hat{I}^-$  or $i=j \in \hat{I}^1$ follow from the above equations noting that $\varepsilon(i,j)=\varepsilon(j,i)$ and, by \eqref{cc1} and \eqref{cc2},  
			\begin{align*}
				&\con{C^\pm_{i}(z)C^\pm_{j}(w)}=(-1)^{\delta_{i,j}}\con{C^\pm_{j}(w)C^\pm_{i}(z)},\\
				&(z-w)\con{C^\pm_{i}(z)C^\mp_{j}(w)}=(-1)^{\delta_{i,j}}(z-w)\con{C^\mp_{j}(w)C^\pm_{i}(z)}.
			\end{align*}
			For example, if $i=m+1$ and $j=m$, let $\epsilon \in\{\pm 1\}$, we have
			\begin{align*}
			    (dqz-w)\con{\Gamma^+_{m+1}(w)\Gamma^+_{m}(z)}\con{C^-_{m}(d^{m-1}q^\epsilon w)C^+_{m}(d^{m}z)} &=d^{-1/2}\dfrac{(dqz-w)(w-dq^{-1}z)}{(d^{m-1}q^\epsilon w-d^mz)}\\
			    &=d^{1/2-m}(dz-q^{-\epsilon}w),
			\end{align*}
			and
			\begin{align*}
			    (dz-qw)\con{\Gamma^+_{m}(z)\Gamma^+_{m+1}(w)}\con{C^+_{m}(d^{m}z)C^-_{m}(d^{m-1}q^\epsilon w)} &=d^{1/2}\dfrac{(dz-qw)(z-d^{-1}q^{-1}w)}{(d^{m}z- d^{m-1}q^\epsilon w)}\\
			    &=d^{1/2-m}(dz-q^{-\epsilon}w).
			\end{align*}
			This shows $(dz-qw)E_m(z)E_{m+1}(w)=(dqz-w)E_{m+1}(w)E_m(z)$.
			
			If $i, j \in \hat{I}^1$ and $i\neq j$, we have $\varepsilon(i,j)=(-1)^{\delta_{m,1}+1}\varepsilon(j,i)$. Thus,
			\begin{align*}
				E_i(z)E_j(w)=(-1)^{\delta_{n,1}+\delta_{m,1}+1}E_j(w)E_i(z) \left(\frac{d^{M_{i,j}}z-q^{-1}w}{w-d^{M_{i,j}}q^{-1}z}\right)^{A_{i,j}}. 
			\end{align*}
			    
			Therefore, the $E$-$E$ relations hold for all $i,j\in \hat{I}$. 
			
			The $F$-$F$ relations are analogous.
			
			The $E$-$F$ relations are trivial for $|i- j|>1$. 
			For $i$ or $j \in \hat{I}^+$ with $i\neq j$, it follows directly from \eqref{gijm}.
			
			If $i \in  \hat{I}^-\cup\{m\}$ and $j=i+ 1$, we have
			\begin{align*}
				  & \con{E_i(z)F_{i+1}(w)}=\con{F_{i+1}(w)E_i(z)}\\ &=d^{4m-2i-1}\con{\Gamma^+_i(z)\Gamma^-_{i+1}(w)}\con{C^+_i(d^{2m-i}z)C^+_i(d^{2m-i-1}w)}=d^{6m-3i-3/2}. 
			\end{align*}
			Thus, $[E_i(z), F_{i+1}(w)]=0$.
			
			The case $i \in  \hat{I}^-\cup\{0\}$ and $j=i-1$ is treated similarly. Due to the presence of the q-difference operators $\partial_z$ and $\partial_w$ in a non-trivial contraction, the expansion of $:E_i(z)F_j(w):$ has four normal-ordered summands. However, the coefficient of each summand is a Laurent polynomial. Thus, $[E_i(z), F_{i-1}(w)]=0$.

			If $i= j \in  \hat{I}^+$, we have
			\begin{align*}
				  & \con{\Gamma^+_i(z)\Gamma^-_i(w)}=\left( \frac{1}{(z-qw)(z-q^{-1}w)} \right)  & (|z|\gg|w|), \\
				  & \con{\Gamma^-_i(w)\Gamma^+_i(z)}= \left( \frac{1}{(z-qw)(z-q^{-1}w)} \right) & (|w|\gg|z|). \\
			\end{align*}
			
			We can change the region where the second rational function is expanded to the same region as the first one by adding $\delta$-functions 
			
			\begin{align}{\label{dchange}}
				  & \left( \frac{1}{(z-qw)(z-q^{-1}w)} \right)=                                                                                                                     &(|w|\gg|z|)  \\
				  & =\left( \frac{1}{(z-qw)(z-q^{-1}w)} \right)- \frac{1}{qw^2(q-q^{-1})} \delta\left(q\frac{w}{z}\right) -\frac{1}{qz^2(q^{-1}-q)} \delta\left(q\frac{z}{w}\right) &(|z|\gg|w|). 
			\end{align}
			
			Now, 
			\begin{align*}
				  & \frac{1}{qw^2(q-q^{-1})} \delta\left(q\frac{w}{z}\right) :\Gamma^+_i(z)\Gamma^-_i(w):                                                             
				=\frac{1}{(q-q^{-1})} \delta\left(q\frac{w}{z}\right)K_i^+(w),\\
				  & \frac{1}{qz^2(q^{-1}-q)} \delta\left(q\frac{z}{w}\right):\Gamma^+_i(z)\Gamma^-_i(w):=-\frac{1}{(q-q^{-1})} \delta\left(q\frac{z}{w}\right)K_i^-(z). 
			\end{align*}
			Therefore, the $E$-$F$ relations follow for $i=j \in \hat{I}^+ $.
			
			For $i=j=0 $  we have
			\begin{align*}
				  & E_0(z)F_0(w)=d^{m-n}:\Gamma^+_0(z)\Gamma^-_0(w):\, \partial_{z}[C^-_{m+n-1}(d^{m-n}z)]C^+_{m+n-1}(d^{m-n}w),  \\
				  & F_0(w)E_0(z)=d^{m-n}:\Gamma^-_0(w)\Gamma^+_0(z):\, C^+_{m+n-1}(d^{m-n}w)\, \partial_{z}[C^-_{m+n-1}(d^{m-n}z)]. 
			\end{align*}
			
			By \eqref{cc2}, we have
			\begin{align*}
			   &d^{m-n} C^-_{m+n-1}(d^{m-n}z)C^+_{m+n-1}(d^{m-n}w)=\dfrac{1}{z-w} &(|z|\gg|w|),\\
			   &d^{m-n} C^+_{m+n-1}(d^{m-n}w)C^-_{m+n-1}(d^{m-n}z)=\dfrac{1}{w-z} &(|w|\gg|z|).
			\end{align*}
			
			Then,
			\begin{align*}
				[E_0(z),F_0(w)]&= :\Gamma^+_0(z)\Gamma^-_0(w):\partial_{z}\left[ \frac{1}{w}\delta\left(\frac{w}{z}\right):C^+_{m+n-1}(d^{m-n}w)C^-_{m+n-1}(d^{m-n}z):\right]. 
			\end{align*}
			
			Now, for all $j\in \hat{I}^-\cup\{m\}$, we have
			\begin{align*}
			    \delta\left(\frac{w}{z}\right):C^+_{j}(w)C^-_{j}(z):=\delta\left(\frac{w}{z}\right)\exp(c^-_{j}(w)-c^-_{j}(z))\exp(c^+_{j}(z)-c^+_{j}(w))=\delta\left(\frac{w}{z}\right).
			\end{align*}
			
			Thus,
			\begin{align*}
				[E_0(z),F_0(w)]&= :\Gamma^+_0(z)\Gamma^-_0(w):\partial_{z}\left[ \frac{1}{w}\delta\left(\frac{w}{z}\right)\right] \\
				&=\frac{1}{zw(q-q^{-1})}\left(\delta\left(q\frac{w}{z}\right)-\delta\left(q\frac{z}{w}\right)\right):\Gamma^+_0(z)\Gamma^-_0(w):. 
			\end{align*}
			Therefore, the $E$-$F$ relations also follow for $i =0 $. The case $i=m$ is analogous and the case $i \in \hat{I}^-$ is longer, but checked by the same procedure.
			
			In any admissible representation, it is enough to check the quadratic relations, then 
			the Serre relations follow automatically. Namely, the Serre relations are checked by commuting each summand and passing to a common region of convergence of the rational functions by adding suitable $\delta$-functions. We check the quartic relation \eqref{Serre3} with $i=m$ as an example.
			
			Write the ten summands in \eqref{Serre3} as follows. Let
			\begin{align}
				E_m(z_1)E_{m+1}(w_1)E_m(z_2)E_{m-1}(w_2)=\mathbf{E}. \label{s1} 
			\end{align} Then, using the $E$-$E$ relations, write the remaining terms of \eqref{Serre3} in the form
			{	\allowdisplaybreaks
			\begin{align}
				& E_{m+1}(w_1)E_m(z_2)E_{m-1}(w_2)E_m(z_1)     =-\frac{(d z_1-q^{-1} w_2)(d z_1-q w_1)}{(d q^{-1}z_1-w_2)(dqz_1-w_1)} \mathbf{E} \label{s2},                                              \\
				& E_m(z_1)E_{m-1}(w_2)E_m(z_2)E_{m+1}(w_1)     =\frac{(d z_2-q^{-1} w_2)(dq z_2- w_1)}{(dz_2-qw_1)(dq^{-1}z_2-w_2)} \mathbf{E},\label{s3}                                                 \\
				& E_{m+1}(w_1)E_m(z_1) E_{m-1}(w_2) E_m(z_2)   =\frac{(d z_2-q^{-1} w_2)(dz_1-qw_1)}{(dq^{-1}z_2-w_2)(dqz_1-w_1)} \mathbf{E},\label{s4}                                                   \\
				& E_{m-1}(w_2) E_m(z_1) E_{m+1}(w_1) E_m(z_2)  =\frac{(d z_2-q^{-1} w_2)(dz_1-q^{-1}w_2)}{(dq^{-1}z_2-w_2)(dq^{-1}z_1-w_2)} \mathbf{E},\label{s5}                                         \\
				& E_{m-1}(w_2) E_m(z_2)E_{m+1}(w_1)E_m(z_1)    =-\frac{(d z_1-q^{-1} w_2)(d z_1-q w_1)(dqz_2-w_1)(d z_2-q^{-1} w_2)}{(d q^{-1}z_1-w_2)(dqz_1-w_1)(dz_2-qw_1)(dq^{-1}z_2-w_2)} \mathbf{E}, \\
				& E_m(z_2)E_{m+1}(w_1)E_m(z_1) E_{m-1}(w_2)    =-\frac{(d z_1-q^{-1} w_2)(d z_1-q w_1)(dqz_2-w_1)}{(dqz_1-w_1)(dz_2-qw_1)(dz_1-q^{-1}w_2)} \mathbf{E},                                    \\
				& E_m(z_2) E_{m-1}(w_2) E_m(z_1)E_{m+1}(w_1)   =-\frac{(d z_1-q^{-1} w_2)(dqz_2-w_1)(dq z_1-w_1)}{(d q^{-1}z_1-w_2)(dqz_1-w_1)(dz_2-qw_1)} \mathbf{E},  \\                       &-[2]\, E_m(z_1)E_{m+1}(w_1)E_{m-1}(w_2)E_m(z_2)  =-[2]\,\frac{(d z_2-q^{-1} w_2)}{(dq^{-1}z_2-w_2)} \mathbf{E},                                                          \\
			&	-[2]\, E_m(z_2)E_{m+1}(w_1)E_{m-1}(w_2)E_m(z_1)  =[2]\,\frac{(d z_1-q^{-1} w_2)(d z_1-q w_1)(dqz_2-w_1)}{(d q^{-1}z_1-w_2)(dqz_1-w_1)(dz_2-qw_1)} \mathbf{E}.\label{s10} 
			\end{align}}
			
			The rational functions of the r.h.s. of the equations \eqref{s2}-\eqref{s10} are expanded in the region given by the increasing order of appearance of the coordinates in the l.h.s.. For example, the rational function in equation \eqref{s2} is expanded in the region  $|w_1|\gg|z_2|\gg|w_2|\gg|z_1|$.
			
			Summing up l.h.s. of equations \eqref{s1}-\eqref{s10} we get the expansion of \eqref{Serre3}. The sum of the rational functions in the r.h.s. vanishes as a rational function. However, similar to $E$-$F$ relation, we must switch to the common convergence region and verify that the coefficients of the delta functions yielded also vanish at the respective support, cf., \eqref{dchange}.
			
			For example, we choose $|z_2|\gg|w_2|\gg|z_1|\gg|w_1|$ as a common region.  The rational function in the r.h.s. \eqref{s4} in this region becomes
			\begin{align*}
				\frac{(d z_2-q^{-1} w_2)(dz_1-qw_1)}{(dq^{-1}z_2-w_2)(dqz_1-w_1)}+(1-q^{-2} )(1-q^2)\delta\left(\frac{dqz_1}{w_1}\right)\delta\left(\frac{dz_2}{qw_2}\right)-                  \\
				-\frac{(d z_2-q^{-1} w_2)(1-q^2)}{q(dq^{-1}z_2-w_2)}\delta\left(\frac{dqz_1}{w_1}\right)-\frac{(1-q^{-2} )(dz_1-qw_1)}{q^{-1}(dqz_1-w_1)}\delta\left(\frac{dz_2}{qw_2}\right). 
			\end{align*}
			Other terms are similar. After changing the region of convergence of all rational functions the $\delta$-functions yielded are
			$\delta\left(\frac{dz_2}{qw_2}\right)$, $\delta\left(\frac{dqz_1}{w_1}\right)$, $\delta\left(\frac{dqz_1}{w_2}\right)$, $\delta\left(\frac{dqz_1}{w_2}\right)\delta\left(\frac{dqz_1}{w_1}\right)$ and $\delta\left(\frac{dz_2}{qw_2}\right)\delta\left(\frac{dqz_1}{w_1}\right)$, and the coefficient of each one vanishes at the respective suport. 
			
			Thus, \eqref{Serre3} with $i=m$ is proved.
		\end{proof}
	\end{thm}

	\bigskip
	\subsection{Screenings}
	
	The $\E_{m|n}$-modules obtained in Theorem \ref{T1} are not irreducible in general. To find their irreducible quotient, we follow \cite{K1} and \cite{KSU}, and introduce the following  $\xi$ - $\eta$ system.
	
	We set $\Res_{z}(\sum_{i\in\Z} a_iz^{-i})=a_1$.
	
	For $i\in \hat{I}^-\cup \{m\}$, introduce the screening operators
	\begin{align*}
		  & \xi_i= \Res_{z}\left(z^{-1}C^-_i(z)\right), \\
		  & \eta_i=\Res_{z}C^+_i(z),                    
	\end{align*}
	acting on $\mathcal{F}_{\Lambda}$, with $\Lambda=\Lambda_j,\, j\not \in \hat{I}^1$, or $\Lambda=(1-a)\Lambda_0+a\Lambda_m$, $a\in \Z$.
	
	The odd operators $\xi_i, \eta_i$, satisfy
	\begin{align*}
		&[\xi_i,\eta_j]=\delta_{i,j},\\
		&[\xi_i,\xi_j]=[\eta_i,\eta_j]=0,\\
		&\mathcal{F}_{\Lambda}=\xi_i\eta_i\mathcal{F}_{\Lambda}\oplus \eta_i\xi_i\mathcal{F}_{\Lambda},
	\end{align*}
	for all $i,j\in \hat{I}^-\cup \{m\}$.
	
	Define 
	\begin{align*}
		\xi=\prod_{i\in \hat{I}^-\cup \{m\}} \xi_i\,,\quad  \eta=\prod_{i\in \hat{I}^-\cup \{m\}} \eta_i. 
	\end{align*}
	
	\begin{prop}
		If $\Lambda=\Lambda_i,\, i\not \in \hat{I}^1$ or $\Lambda=(1-a)\Lambda_0+a\Lambda_m, a\in \Z$, the  screening operators $\eta_i, i\in \hat{I}^-\cup \{m\},$ (super)commute with the $\E_{m|n}$-action on $\mathcal{F}_{\Lambda}$ given by Theorem \ref{T1}. Thus, $\ker{\eta}$ and $\coker{\eta}$ are $\E_{m|n}$-modules.
		
		\begin{proof}
			It is enough to show $[\eta_i,C^+_i(w)]=[\eta_i,\partial_{w}C^-_i(w)]=0$.
			
			Using \eqref{cc1} we have
			\begin{align*}
				[\eta_i,C^+_i(w)]=\Res_{z}[C^+_i(z),C^+_i(w)]=0, 
			\end{align*}
			and by \eqref{cc2}
			\begin{align*}
				[\eta_i,\partial_{w}C^-_i(w)]=\partial_{w}\left(\Res_{z}[C^+_i(z),C^-_i(w)]\right)=\partial_{w}\left(1\right)=0. 
			\end{align*}
		\end{proof}
	\end{prop}
	
	Level $1$ partially integrable representations of $U_q\,\widehat{\mathfrak{sl}}_{m|n}$ with $m\neq n$ were constructed in  \cite{KSU} using the formulas in Theorem \ref{T1} for $i\in I$ and $d=1$. Our space $\mathcal{F}_{\Lambda}$ differs from theirs by the extra current $H_0(z)$ present in $U^{ver}_q\,\widehat{\mathfrak{gl}}_{m|n}$.  The conjectural identification given in \cite{K2} and \cite{KSU} in our context is the following.
	
	\begin{conj}{\label{conjM}} We have the following identifications
		\begin{align*}
			  & V(\Lambda_i)=\ker\eta=\eta\xi\mathcal{F}_{\Lambda_i}                                   & (i\in I),                      \\
			  & V((1-a)\Lambda_0+a\Lambda_m)=\mathcal{F}_{(1-a)\Lambda_0+a\Lambda_m}                   & (a\in \mathbb{C}\setminus \Z), \\
			  & V((1-a)\Lambda_0+a\Lambda_m)=\coker\eta=\xi\eta\mathcal{F}_{(1-a)\Lambda_0+a\Lambda_m} & (a\in \Z_{>0}),                \\
			  & V((1-a)\Lambda_0+a\Lambda_m)=\ker\eta=\eta\xi\mathcal{F}_{(1-a)\Lambda_0+a\Lambda_m}   & (a\in \Z_{\leq 0}),            
		\end{align*}
		where $V(\Lambda)$ is the irreducible highest weight $U^{ver}_q\widehat{\mathfrak{gl}}_{m|n}$-module with highest  weight $\Lambda$. 
	\end{conj}

	\section{Evaluation homomorphism}{\label{ev}}
	
	In this section, we construct an evaluation map from $\E_{m|n}$ to a suitable completion $\Uge$ of $\Ug$. See Appendix \ref{Asl}.  We follow the  strategy used in \cite{FJM2}. 
	
	\subsection{Fused Currents}
	
	Introduce the following fused currents in $\Uge$
	\begin{align*}
		\textsf{X}^+(z)=  & \left[\prod_{i=1}^{m+n-2}\left(1-\dfrac{z_i}{z_{i+1}}\right)\right]x^+_{m+n-1}(q^{-m+n-1}c^{-1}z_{m+n-1})\cdots x^+_{m+i}(q^{-m+i}c^{-1}z_{m+i})\cdots \times  \\
		                  & \times x^+_{m}(q^{-m}c^{-1}z_{m}) \cdots x^+_{i}(q^{-i}c^{-1}z_{i})\cdots x^+_{1}(q^{-1}c^{-1}z_{1})\Bigl|_{z_1=\cdots=z_{m+n-1}=z}\;,                                                                              \\
		\textsf{X}^-(z)=  & \left[\prod_{i=1}^{m+n-2}\left(1-\dfrac{z_{i+1}}{z_{i}}\right)\right]x^-_{1}(q^{-1}c^{-1}z_{1})\cdots x^-_{i}(q^{-i}c^{-1}z_{i}) \cdots x^-_{m}(q^{-m}c^{-1}z_{m})\times                \\
		                  & \times\cdots x^-_{m+i}(q^{-m+i}c^{-1}z_{m+i})\cdots x^-_{m+n-1}(q^{-m+n-1}c^{-1}z_{m+n-1})\Bigl|_{z_1=\cdots=z_{m+n-1}=z}\;,                                                            \\
		\textsf{k}^\pm(z) & =\prod_{i=1}^{m}k_i^\pm(q^{-i} c^{-1}z)\prod_{j=m+1}^{m+n-1}k_{j}^\pm(q^{-2m+j} c^{-1}z).                                                                                               
	\end{align*}
	See \cite{FJMM2} for the details on fused currents.
	
	The homomorphism $v$ defined in the Lemma \ref{lemv} maps the element $h_{0,r}$ in the following way
	\begin{align*}
		v(h_{0,r})=\frac{1}{\beta_{0,r}}\left(\gamma_{0,r}H_{0,r} + \sum_{i\in \hat{I}^+\cup \{m\}}(\gamma_{i,r}-\beta_{i,r}d^{ir})H_{i,r} + \sum_{j\in  \hat{I}^-}(\gamma_{j,r}-\beta_{j,r}d^{(2m-j)r})H_{j,r}\right), 
	\end{align*}
	where $\{\gamma_{i,r}\}_{i\in \hat{I}}$ and $\{\beta_{i,r}\}_{i\in \hat{I}}$ are the fixed solutions of the systems \eqref{sys1} and \eqref{sys2}, respectively. 
	
	For each $r\in \Zz$, define
	\begin{align*}
		\tilde{h}_{0,r}=\frac{1}{\gamma_{0,r}}\left( \beta_{0,r} h_{0,r} + \sum_{i\in \hat{I}^+\cup \{m\}}(\beta_{i,r}-\gamma_{i,r}d^{-ir})h_{i,r} + \sum_{j\in  \hat{I}^-}(\beta_{j,r}-\gamma_{j,r}d^{-(2m-j)r})h_{j,r}\right). 
	\end{align*}
	Thus, $ v(\tilde{h}_{0,r})=H_{0,r}$ for all $r\in \Zz$.
	
	Define $A_{\pm r}, B_{\pm r}, r\in \Z_{>0} $ by 
	\begin{align*}
		  & A_r=-\frac{q-q^{-1}}{c^r-c^{-r}}\left(\Tilde{h}_{0,r}+\sum_{i=1}^m (c^2q^{i})^{r}h_{i,r}+\sum_{j=m+1}^{m+n-1} (c^2q^{2m-j})^{r}h_{j,r}\right),               \\
		  & A_{-r}=\frac{q-q^{-1}}{c^r-c^{-r}}\;c^{-r}\left(\Tilde{h}_{0,-r}+\sum_{i=1}^m q^{-ir}h_{i,-r}+\sum_{j=m+1}^{m+n-1} q^{(-2m+j)r}h_{j,-r}\right),              \\
		  & B_{r}=\frac{q-q^{-1}}{c^r-c^{-r}}\;c^{r}\left(\Tilde{h}_{0,r}+\sum_{i=1}^m q^{ir}h_{i,r}+\sum_{j=m+1}^{m+n-1} q^{(2m-j)r}h_{j,r}\right),                     \\
		  & B_{-r}=-\frac{q-q^{-1}}{c^r-c^{-r}}\left(\Tilde{h}_{0,-r}+\sum_{i=1}^m (c^{-2}q^{-i})^{r}h_{i,-r}+\sum_{j=m+1}^{m+n-1} (c^{-2}q^{-2m+j})^{r}h_{j,-r}\right), 
	\end{align*}
	and let $A^{\pm}(z)=\sum_{r>0}A_{\pm r}z^{\mp r}$, $B^{\pm}(z)=\sum_{r>0}B_{\pm r}z^{\mp r}$.
	
	Let also $\mathcal{K}=q^{-\Lambda_{m+n-1}-\Lambda_{1}}$. We have $\mathcal{K}x_i^\pm(z)\mathcal{K}^{-1}=q^{\mp(\delta_{1,i}+\delta_{m+n-1,i})}x_i^\pm(z)$.
	\begin{thm}{\label{T2}}
		Fix $u\in \mathbb{C}^\times$. The following map is a surjective homomorphism of superalgebras $\eva_u : \E_{m|n} \rightarrow \Uge$ with  $ C^2=q_3^{m-n}$:
		\begin{align*}
			&K\mapsto 1,\quad C\mapsto c,\quad H^{ver}(z)\mapsto h(z),\\
			  & E_i(z)\mapsto x^+_i(d^iz),\quad F_i(z)\mapsto x^-_i(d^iz),\quad K_i^\pm(z)\mapsto k_i^\pm(d^iz) \quad                & (i\in \hat{I}^+\cup\{m\}), \\ 
			  & E_j(z)\mapsto x^+_j(d^{2m-j}z),\quad F_j(z)\mapsto x^-_j(d^{2m-j}z),\quad K_j^\pm(z)\mapsto k_j^\pm(d^{2m-j}z) \quad & (j\in \hat{I}^-),          \\ 
			&E_0(z)\mapsto u^{-1}e^{A_{-}(z)}\textsf{X}^-(z)e^{A_{+}(z)}\mathcal{K},\\
			&F_0(z)\mapsto u\,\mathcal{K}^{-1}e^{B_{-}(z)}\textsf{X}^+(z)e^{B_{+}(z)}.
		\end{align*}
		
		Moreover, the evaluation map $\eva_u$ is graded: if $X\in\E_{m|n}$ and $\deg(X)=\left(d_0,d_1,\dots,d_{m+n-1};d_{\delta}\right)$, then $\deg(ev_u(X))=\left(d_1-d_0,\dots,d_{m+n-1}-d_0;d_{\delta}\right)$.
		
		\begin{proof}
			For simplicity, we fix $u=1$ and write $\eva_1=\eva$.
			The relations with no  index $0$ are clear.

			For $i \in \hat{I}, r>0$, we have 
			\begin{align*}
				  & [h_{i,r},e^{A^+(z)}]=0,                                                                                                                                                           \\
				  & [h_{i,r},e^{A^-(z)}]=z^re^{A^-(z)}c^{-r}\left( \sum_{j\in \hat{I}^+\cup \hat{I}^1}\frac{[rA_{i,j}]}{r}q^{-jr}+\sum_{j\in \hat{I}^-}\frac{[rA_{i,j}]}{r}q^{-(2m-j)r}\right),       \\
				  & [h_{i,r},\textsf{X}^-(z)]=-z^r\textsf{X}^-(z)c^{-r}\left( \sum_{j\in \hat{I}^+\cup\{m\}}\frac{[rA_{i,j}]}{r}q^{-jr}+\sum_{j\in \hat{I}^-}\frac{[rA_{i,j}]}{r}q^{-(2m-j)r}\right). 
			\end{align*}
			Thus,
			\begin{align*}
				  & \eva([H_{i,r},E_0(z)])=z^rc^{-r}\frac{[rA_{i,0}]}{r}\eva(E_0(z)). 
			\end{align*}
			The $H$-$E$ relations with $r<0$ and the $H$-$F$ relations can be checked in the same way.
			
			To check the $E$-$E$ relations we first use \eqref{cont1} and \eqref{cont2} to get
			\begin{align*}
				  & \eva(E_0(z)E_i(w))=e^{A^{-}(z)}\textsf{X}^-(z)\eva(E_i(w))e^{A^{+}(z)}\mathcal{K}\left(\frac{z-q_3^{-1}w}{z-q_1w}\right)^{\delta_{i,1}}q^{-\delta_{i,m+n-1}-\delta_{i,1}}, \\
				  & \eva(E_i(w)E_0(z))=e^{A^{-}(z)}\eva(E_i(w))\textsf{X}^-(z)e^{A^{+}(z)}\mathcal{K}\left(\frac{w-q_3z}{w-q_1^{-1}z}\right)^{\delta_{i,m+n-1}}.                               
			\end{align*}
			
			For $i\neq 1, m+n-1$, $\eva([E_0(z),E_i(w)])=0$ by \eqref{f1}.
			
			For $i=1$, the $E$-$E$ relation reduces to 
			\begin{align*}
				\left(q^{-1}z-dw\right) e^{A^{-}(z)}[\textsf{X}^-(z)x^+_1(dw)]e^{A^{+}(z)}\mathcal{K}=0, 
			\end{align*}
			which follows from \eqref{f6}. The case $i=m+n-1$ is similar. 
			The case $i=0$ is checked using \eqref{cont3}, \eqref{cont4} and \eqref{cont9}.
			
			The $F$-$F$ relations are verified by the same argument.
			
			For the $E$-$F$ relations
			$$\eva([E_0(z),F_i(w)])=0\quad (i\neq 0),$$
			we proceed as in the $E$-$E$ case by bringing $A^-(z)$ to the left and $A^+(z)$ to the right using \eqref{cont3} and \eqref{cont4}. The relations then follow from \eqref{f1}, \eqref{f4} and \eqref{f5}. The same is done for $\eva([E_i(z),F_0(w)])=0\quad (i\neq 0)$.
			
			For the $i=0$ case, using \eqref{cont1},\eqref{cont6} and \eqref{cont10}, we get
			\begin{align*}
				\eva(E_0(z)F_0(w))=e^{A^{-}(z)}e^{B^{-}(w)}\textsf{X}^-(z)\textsf{X}^+(w)e^{A^{+}(z)}e^{B^{+}(w)}, 
			\end{align*}
			and similarly
			\begin{align*}
				\eva(F_0(w)E_0(z))=e^{A^{-}(z)}e^{B^{-}(w)}\textsf{X}^+(w)\textsf{X}^-(z)e^{A^{+}(z)}e^{B^{+}(w)}. 
			\end{align*}
			By \eqref{f8}, 
			\begin{align*}
				\eva([E_0(z),F_0(w)])=e^{A^{-}(z)}e^{B^{-}(w)}\frac{1}{q-q^{-1}}\Bigl(\delta\left(c\frac{w}{z}\right)\textsf{k}^-(w)-\delta\left(c\frac{z}{w}\right)\textsf{k}^+(z) \Bigr)e^{A^{+}(z)}e^{B^{+}(w)}. 
			\end{align*}
			
			The relation \eqref{EF} with $i=j=0$ follows from
			\begin{align*}
				e^{A^{-}(z)}e^{B^{-}(cz)}=\Tilde{k}_0^-(z),\quad   e^{A^{-}(cw)}e^{B^{-}(w)}=k_0(\textsf{k}^-(w))^{-1}, \\
				e^{A^+(cw)}e^{B^+(w)}=\Tilde{k}_0^+(w),\quad e^{A^+(z)}e^{B^+(cz)}=k_0^{-1}(\textsf{k}^+(z))^{-1},      
			\end{align*}
			where $\Tilde{k}_0^\pm=\exp\Bigl(\pm(q-q^{-1})\sum_{r>0}\Tilde{h}_{0,\pm r}z^{\mp r}\Bigr)$.
			
			Finally, we check the Serre relations.
			
			For the relation
			\begin{align*}
			    \eva\left(\Sym_{{z_1,z_2}}[E_1(z_1),[E_1(z_2),E_{0}(w)]_q]_{q^{-1}}\right)=0
			\end{align*}
			we use \eqref{cont1} and \eqref{f6}  to obtain
			\begin{align*}
				  &-(q+q^{-1}) \eva\left(E_1(z_1)E_{0}(w)E_1(z_2)\right)= -q(q+q^{-1})\left(\frac{z_2-q_3 w}{z_2-q_1^{-1}w}\right)\eva\left(E_1(z_1)E_1(z_2)E_{0}(w)\right),                                                 \\
				  & \eva\left(E_{0}(w)E_1(z_1)E_1(z_2)\right)=q^{2}\left(\frac{z_1-q_3 w}{z_1-q_1^{-1}w}\right)\left(\frac{z_2-q_3 w}{z_2-q_1^{-1}w}\right) \eva\left(E_1(z_1)E_1(z_2)E_{0}(w)\right). 
			\end{align*}
			
			Thus, 
			\begin{align*}
				 &\Sym_{{z_1,z_2}} (\eva\left[E_1(z_1),[E_1(z_2),E_{0}(w)]_q]_{q^{-1}}\right)=\\
				 &\frac{(1-q^2)w}{q_3(w-q_1 z_1)(w-q_1 z_2)}\,\Sym_{{z_1,z_2}} \left((z_1-q^2z_2)\eva\left(E_1(z_1)E_1(z_2)E_{0}(w)\right)\right)=0,
			\end{align*}
			where the last equality follows from the quadratic relation for $x^+_1(dz_1)x^+_1(dz_2)$.
			
			The Serre relations in all the remaining cases are checked in the same way. 
			
			The statement about grading is straightforward.
		\end{proof}
	\end{thm}
	
	By Theorem \ref{T2}, any admissible $U_q\,\widehat{\mathfrak{gl}}_{m|n}$-module of generic level $c$ can be pulled back by $\eva_u$ to a representation of $\E_{m|n}$ with $q_3^{m-n}=c^2$ and $q_2=q^2$.
Such $\E_{m|n}$ modules are called {\it evaluation modules}.
	
	\medskip
There exist another evaluation map $\eva_u^{(1)}$ obtained by composing the map $\eva_u$ for $\E_{n|m}$ with the change of parity isomorphism \eqref{cph}. For this map we have $C^2=q_1^{m-n}$.

\medskip
	
The evaluation maps $\eva_u$ and $\eva_u^{(1)}$ prefer our choice of the vertical subalgebra related to the zero node of the Dynkin diagram. In fact, there are evaluation maps which prefer any node. The ones related to odd node $m$ are obtained by the diagram automorphism \eqref{diagiso}. However, the vertical subalgebras related to even nodes appear in non-standard parities, and we do not discuss the corresponding isomorphisms or evaluation maps in this paper.

	\medskip

	\appendix
	\allowdisplaybreaks
	
	\section{}\label{app1}
	
	In this Appendix, we collect some useful formulas for commutation relations of various currents.

	\begin{lem} For $i\in I$, we have
		\begin{align}
			  & e^{A^+(z)}x^+_i(dw)e^{-A^+(z)}=x^+_i(dw)\left(\frac{z-q_3^{-1}w}{z-q_1w}\right)^{\delta_{i,1}},\label{cont1}                                 \\
			  & e^{-A^-(z)}x^+_i(d^{m-n+1}w)e^{A^-(z)}=x^+_i(d^{m-n+1}w)\left(\frac{w-q_3z}{w-q_1^{-1}z}\right)^{\delta_{i,m+n-1}},\label{cont2}             \\
			  & e^{A^+(z)}x^-_i(dw)e^{-A^+(z)}=x^-_i(dw)\left(\frac{z-cq_1w}{z-cq_3^{-1}w}\right)^{\delta_{i,1}},\label{cont3}                               \\
			  & e^{-A^-(z)}x^-_i(d^{m-n+1}w)e^{A^-(z)}=x^-_i(d^{m-n+1}w)\left(\frac{w-cq_1^{-1}z}{w-cq_3z}\right)^{\delta_{i,m+n-1}},\label{cont4}           \\
			  & e^{B^+(z)}x^+_i(d^{m-n+1}w)e^{-B^+(z)}=x^+_i(d^{m-n+1}w)\left(\frac{z-c^{-1}q_1w}{z-c^{-1}q_3^{-1}w}\right)^{\delta_{i,m+n-1}},\label{cont5} \\
			  & e^{-B^-(z)}x^+_i(dw)e^{B^-(z)}=x^+_i(dw)\left(\frac{w-c^{-1}q_1^{-1}z}{w-c^{-1}q_3z}\right)^{\delta_{i,1}},\label{cont6}                     \\
			  & e^{B^+(z)}x^-_i(d^{m-n+1}w)e^{-B^+(z)}=x^-_i(d^{m-n+1}w)\left(\frac{z-q_3^{-1}w}{z-q_1w}\right)^{\delta_{i,m+n-1}},\label{cont7}             \\
			  & e^{-B^-(z)}x^-_i(dw)e^{B^-(z)}=x^-_i(dw)\left(\frac{w-q_3z}{w-q_1^{-1}z}\right)^{\delta_{i,1}},\label{cont8}                                 \\
			  & e^{A^+(z)}e^{A^-(w)}=e^{A^-(w)}e^{A^+(z)}\frac{(z-w)^2}{(z-q_2w)(z-q_2^{-1}w)},\label{cont9}                                                 \\
			  & e^{A^+(z)}e^{B^-(w)}=e^{B^-(w)}e^{A^+(z)}\frac{(z-cq_2w)(z-c^{-1}q_2^{-1}w)}{(z-cw)(z-c^{-1}w)},\label{cont10}                               \\
			  & e^{B^+(z)}e^{B^-(w)}=e^{B^-(w)}e^{B^+(z)}\frac{(z-w)^2}{(z-q_2w)(z-q_2^{-1}w)},\label{cont11}                                                \\
			  & e^{B^+(w)}e^{A^-(z)}=e^{A^-(z)}e^{B^+(w)}\frac{(z-c^{-1}q_2w)(z-cq_2^{-1}w)}{(z-cw)(z-c^{-1}w)}.\label{cont12}                               
		\end{align}
		\qed    
	\end{lem}
	
	\begin{lem} The fused currents $\textsf{X}^{\pm}(z)$ satisfy
		\begin{align}
			  & [x^{\pm}_i(z),\textsf{X}^{\pm}(w)]=[x^{\pm}_i(z),\textsf{X}^{\mp}(w)]=0 \qquad (i\neq 1,m+n-1),\label{f1}                                                                      \\
			  & q(w-c^{-1}q_3z)\textsf{X}^{+}(z)x^{+}_1(dw)=(w-c^{-1}q_1^{-1}z)x^{+}_1(dw)\textsf{X}^{+}(z),\label{f2}                                                                         \\  
			  & q(z-c^{-1}q_1w)\textsf{X}^{+}(z)x^{+}_{m+n-1}(d^{m-n+1}w)=(z-c^{-1}q_3^{-1}w)x^{+}_{m+n-1}(d^{m-n+1}w)\textsf{X}^{+}(z),\label{f3}                                             \\  
			  & q(z-cq_1w)\textsf{X}^{-}(z)x^{-}_1(dw)=(z-cq_3^{-1}w)x^{-}_1(dw)\textsf{X}^{-}(z),\label{f4}                                                                                   \\  
			  & q(w-cq_3z)\textsf{X}^{-}(z)x^{-}_{m+n-1}(d^{m-n+1}w)=(w-cq_1^{-1}z)x^{-}_{m+n-1}(d^{m-n+1}w)\textsf{X}^{-}(z),\label{f5}                                                       \\  
			  & (z-q_3^{-1}w)[\textsf{X}^{-}(z),x^{+}_1(dw)]=(w-q_3z)[\textsf{X}^{-}(z),x^{+}_{m+n-1}(d^{m-n+1}w)]=0,\label{f6}                                                                \\  
			  & (w-q_3z)[\textsf{X}^{+}(z),x^{-}_1(dw)]=(z-q_3^{-1}w)[\textsf{X}^{+}(z),x^{-}_{m+n-1}(d^{m-n+1}w)]=0,\label{f7}                                                                \\  
			  & [\textsf{X}^{+}(w),\textsf{X}^{-}(z)]=\frac{1}{q-q^{-1}}\Bigl(-\delta\left(c\frac{z}{w}\right)\textsf{k}^+(z)+\delta\left(c\frac{w}{z}\right)\textsf{k}^-(w) \Bigr),\label{f8} \\
			  & [\textsf{X}^{\pm}(z),\textsf{X}^{\pm}(w)]=0.\label{f9}                                                                                                                         
		\end{align}
		\qed
	\end{lem}

	\section{}\label{Asl}
	In this Appendix, we set up notation involving the algebra $U_q\,\widehat{\mathfrak{sl}}_{m|n}$.

	The presentations of the superalgebra $U_q\,\widehat{\mathfrak{sl}}_{m|n}$ in Drinfeld-Jimbo and new Drinfeld forms were given in \cite{Y}. We recall them here using the standard choice of parity.
	
	Let $\bar{\alpha}_i$, $\bar{\Lambda}_i$, $i\in I$, be the simple roots and fundamental weights of $\sln_{m|n}$, respectively. Let $\braket{\cdot}{\cdot}$ be the symmetric bilinear form given by $\braket{\bar{\alpha}_{i}}{\bar{\alpha}_{j}}={A}_{i,j}$ and $\braket{\bar{\alpha}_{i}}{\bar{\Lambda}_{j}}=\delta_{i,j}$, $i,j \in I$. Let $\Lambda_0$ be the affine fundamental weight and $\delta$ be the null root of $\widehat{\sln}_{m|n}$. We have  $\braket{\Lambda_0}{\Lambda_0}=\braket{\delta}{\delta}=\braket{\delta}{\bar{\alpha}_i}=\braket{\delta}{\bar{\Lambda}_i}=0$, $i\in I$, and $\braket{\Lambda_0}{\delta}=1$. The remaining affine fundamental weights are $\Lambda_i=\bar{\Lambda}_i+\Lambda_0$, $i\in I$. The affine roots are $\alpha_i=\bar{\alpha}_i$, $i\in I$, and $\alpha_0=\delta-\sum_{i\in I}\alpha_i$. Set also $\bar{\alpha}_0=-\sum_{i\in I}\bar{\alpha}_i$.
	
	In the Drinfeld-Jimbo realization, the algebra $U_q\,\widehat{\mathfrak{sl}}_{m|n}$ is generated by Chevalley type elements $e_i, f_i, t_i,\, i \in \hat{I} $. The parity of generators is given by $|e_0|=|f_0|=|e_m|=|f_m|=1$ and $0$ otherwise. The relations are as follows.
	\begin{align*}
		&t_it_j=t_jt_i,\\
		&t_ie_jt_i^{-1}=q^{A_{i,j}}e_j,\\
		&t_if_jt_i^{-1}=q^{-A_{i,j}}f_j,\\
		&[e_i,f_j]=\delta_{i,j}\frac{t_i-t_i^{-1}}{q-q^{-1}},\\
		  & [e_i,e_j]=[f_i,f_j]=0                                                                 & (A_{i,j}=0),                \\
		  & [e_i,[e_i,e_{i\pm 1}]_{q}]_{q^{-1}}=[f_i,[f_i,f_{i\pm 1}]_{q}]_{q^{-1}}=0             & (i\not\in \hat{I}^1),       \\
		  & [e_i,[e_{i+1},[e_i,e_{i-1}]_q]_{q^{-1}}]=[f_i,[f_{i+1},[f_i,f_{i-1}]_q]_{q^{-1}}]=0   & (i\in \hat{I}^1, mn\neq 2), \\
		  & [e_2,[e_0,[e_{2},[e_0,e_{1}]_q]]]_{q^{-1}}=[e_0,[e_2,[e_{0},[e_2,e_{1}]_q]]]_{q^{-1}} & ((m,n)=(2,1)),              \\
		  & [f_2,[f_0,[f_{2},[f_0,f_{1}]_q]]]_{q^{-1}}=[f_0,[f_2,[f_{0},[f_2,f_{1}]_q]]]_{q^{-1}} & ((m,n)=(2,1)),              \\
		  & [e_1,[e_0,[e_{1},[e_0,e_{2}]_q]]]_{q^{-1}}=[e_0,[e_1,[e_{0},[e_1,e_{2}]_q]]]_{q^{-1}} & ((m,n)=(1,2)),              \\
		  & [f_1,[f_0,[f_{1},[f_0,f_{2}]_q]]]_{q^{-1}}=[f_0,[f_1,[f_{0},[f_1,f_{2}]_q]]]_{q^{-1}} & ((m,n)=(1,2)).              
	\end{align*}
	
	The element $t_0t_1\cdots t_{m+n-1}$ is central.
	
	\medskip
	
	In the new Drinfeld realization, the algebra  $U_q\,\widehat{ \mathfrak{sl} }_{m|n}$ is generated by current generators $x^\pm_{i,n}, h_{i,r}$, $k^{\pm 1}_i, c^{\pm 1}$, $i \in I,\; n\in \Z,\; r\in \Zz$, satisfying
	\begin{align*}
		&\text{$c$ is central},\quad k_ik_j=k_jk_i,\quad k_ix^\pm_j(z)k_i^{-1}=q^{\pm A_{i,j}}x^\pm_j(z),\\
		&[h_{i,r},h_{j,s}]=\delta_{r+s,0}\frac{[rA_{i,j}]}{r}\frac{c^r-c^{-r}}{q-q^{-1}},\\
		&[h_{i,r},x^{\pm}_j(z)]=\pm\frac{[rA_{i,j}]}{r}c^{-(r\pm|r|)/2}z^rx^\pm_j(z),\\
		&[x^+_i(z),x^-_j(w)]=\frac{\delta_{i,j}}{q-q^{-1}}\Bigl(\delta\bigl(c\frac{w}{z}\bigr)k_i^+(w)-\delta\bigl(c\frac{z}{w}\bigr)k_i^-(z)\Bigr),\\
		  & (z-q^{\pm A_{i,j}}w)x^\pm_i(z)x^\pm_j(w)+(w-q^{\pm A_{i,j}}z)x^\pm_j(w)x^\pm_i(z)=0            & (A_{i,j}\neq 0)\,, \\
		  & [x^\pm_i(z),x^\pm_j(w)]=0                                                                      & (A_{i,j}=0)\,,     \\
		  & \Sym_{z_1,z_2}[x^\pm_i(z_1),[x^\pm_i(z_2),x^\pm_{i\pm 1}(w)]_{q}]_{q^{-1}}=0\,                        & (i\neq m)\,,       \\
		  & \Sym_{{z_1,z_2}}[x^\pm_m(z_1),[x^\pm_{m+1}(w_1),[x^\pm_m(z_2),x^\pm_{m-1}(w_2)]_q]_{q^{-1}}]=0 & (m,n>1)\,,         
	\end{align*}
	where $x^\pm_i(z)=\sum_{k\in \Z}x^\pm_{i,k}z^{-k}\,,$ $k_i^\pm(z)=k_i^{\pm 1}\exp \left(\pm (q-q^{-1})\sum_{r>0}h_{i,\pm r}z^{\mp r}\right)$.
	
	An isomorphism between the two realizations is given by
	\begin{align*}
		  & e_i\mapsto x^+_{i,0},\quad f_i\mapsto x^-_{i,0},\quad t_i\mapsto k_i & (i \in I), \\
		&t_0\mapsto c(k_1k_2\cdots k_{m+n-1})^{-1},\\
		&e_0\mapsto (-1)^n(k_1k_2\cdots k_{m+n-1})^{-1}[x^-_{m+n-1,0},\cdots,[x^-_{m+1,0},[x^-_{m,0},\cdots,[x^-_{2,0},x^-_{1,1}]_{q^{-1}}\cdots]_{q^{-1}}]_q\cdots]_q,\\
		&f_0=k_1k_2\cdots k_{m+n-1}[\cdots [[\cdots[x^+_{1,-1},x^+_{2,0}]_q,\cdots x^+_{m,0}]_q,x^+_{m+1,0}]_{q^{-1}},\cdots x^+_{m+n-1,0}]_{q^{-1}}.
	\end{align*}
	Note that $t_0t_1\cdots t_{m+n-1}\mapsto c$.
	
	The quantum affine superalgebra $U_q\,\widehat{\mathfrak{gl}}_{m|n}$ is obtained from $U_q\,\widehat{\mathfrak{sl}}_{m|n}$ in the new Drinfeld realization by including the currents $k_0^{\pm}(z)$ subject to the same relations.
	
	For $r\in \Zz$, let
	\begin{align}{\label{sys2}}
		\beta_{i,r}=\begin{cases}\dfrac{q^{(m-n-i)r}+q^{ir}}{q^r-q^{-r}} & (i\in \hat{I}^+\cup\hat{I}^1),\vspace{.1cm} \\
		\dfrac{q^{(i-m-n)r}+q^{(2m-i)r}}{q^r-q^{-r}}                    & (i\in \hat{I}^-).                       \\
		\end{cases}
	\end{align}
	
	The coefficients $\beta_{i,r}$ are solutions of the system
	\begin{align*}
	    	\sum_{i\in \hat{I}}\beta_{i,r}[rA_{i,j}]=0 \quad (j\in I)\,.
	\end{align*}
	
	Then, the elements $h_r= \sum_{i\in \hat{I}}\beta_{i,r}h_{i,r} \in \Ug$ commute with $\Us \subset \Ug$ and satisfy
	
	$$[h_r,h_s]=\delta_{r+s,0}[(n-m)r]\frac{1}{r}\frac{c^r-c^{-r}}{q-q^{-1}}.$$
	Set $h(z)=\sum_{k\in \Zz}h_kz^{-k}$.
	
	\medskip
	
	We use a completion of $U_q\,\widehat{\mathfrak{gl}}_{m|n}$, denoted by $\widetilde{U}_q\,\widehat{\mathfrak{gl}}_{m|n}$, obtained by performing the following two steps.
	
	Let $\widehat{Q}_{m|n}$ be the root lattice of $\widehat{\mathfrak{sl}}_{m|n}$.
	The algebra $U_q\,\widehat{\mathfrak{gl}}_{m|n}$ contains the group algebra $\C [\widehat{Q}_{m|n}]$ of the root lattice $\widehat{Q}_{m|n}$ generated by  $k_i=q^{\alpha_i}, i \in \hat{I}$. As a first step, we extend it to the weight lattice in a straightforward way. Namely, let $P$ be the $\widehat{\mathfrak{sl}}_{m|n}$ weight lattice and $\C [P]$ the corresponding group algebra spanned by $q^\La, \La\in P$. We have an inclusion of algebras $\C [\widehat{Q}_{m|n}]\subset \C [P]$. Let $U_P$ be the superalgebra $U_{P}=U_q\,\widehat{\mathfrak{gl}}_{m|n}\otimes_{\C [\widehat{Q}_{m|n}]} \C [P]$ with the relations
	\begin{align*}
		q^{\Lambda}q^{\Lambda'}=q^{\Lambda+\Lambda'}\,,\quad                                		q^0=1\,,\quad q^\Lambda x^\pm_i(z)q^{-\Lambda}=q^{\pm\braket{\Lambda}{\alpha_i}}x^\pm_i(z) \quad (\Lambda, \Lambda'\in P). 
	\end{align*}
	
		For each $i\in I$, the superalgebra $U_{P}$ has a $\Z$-grading given by 
	\begin{align*}
		&\deg_i(x^\pm_{j,k})=\pm \delta_{i,j},\quad &\deg_i(h_{j,r})=\deg_i  (q^\Lambda)=\deg_i(c)=0\quad  (j\in I,  k\in \Z, r\in \Zz, \Lambda \in P).
	\end{align*}

	There is also the \textit{homogeneous $\Z$-grading} given by
	\begin{align*}
		\deg_{\delta}( x^\pm_{j,k})=k,\quad \deg_{\delta}(h_{j,r})=r, \quad \deg_{\delta}(q^\Lambda)=\deg_{\delta}(c)=0 \quad  (j\in I,  k\in \Z, r\in \Zz, \Lambda \in P). 
	\end{align*}
	
	Thus the superalgebra $U_P$ has a $\Z^{m+n}$-grading given by
	\begin{align*}
	    \deg (X)=\left(\deg_1(X),\dots,\deg_{m+n-1};\deg_{\delta}(X) \right), \qquad X\in U_P.
	\end{align*}
	
	As the second step, we define $\widetilde{U}_q\,\widehat{\mathfrak{gl}}_{m|n}$ to be the completion of $U_{P}$ with respect to the homogeneous grading in the positive direction. The elements of $\widetilde{U}_q\,\widehat{\mathfrak{gl}}_{m|n}$ are series of the form $\sum_{j=s}^\infty g_j$, with $g_j\in U_P,$ $\deg_{\delta} g_j=j.$
	
	\begin{lem}{\label{lem1}}
		We have an embedding
		\begin{align*}
			U_q\,\widehat{\mathfrak{gl}}_{m|n} \rightarrow \widetilde{U}_q\,\widehat{\mathfrak{gl}}_{m|n}. 
		\end{align*}
		\qed
	\end{lem}
	
	 A $U_q\,\widehat{\mathfrak{gl}}_{m|n}$-module $V$ is {\it admissible} if for any $v\in V$ there exist $N=N_v>0$ such that $xv=0$ for all $x \in U_q\,\widehat{\mathfrak{gl}}_{m|n}$ with $\deg_\delta x> N $. Any admissible $U_q\,\widehat{\mathfrak{gl}}_{m|n}$-module is also an $\widetilde{U}_q\,\widehat{\mathfrak{gl}}_{m|n}$-module. 
	 
	 A $U_q\,\widehat{\mathfrak{gl}}_{m|n}$-module $V$ is called {\it highest weight module} if $V$ is generated by the highest weight vector $v$:
	 $$
	 V=U_q\,\widehat{\mathfrak{gl}}_{m|n}v, \qquad e_iv=0,\qquad k_0^+(z)v=q^{\lambda_0} v,\qquad t_jv=q^{\lambda_j}v, \qquad i\in \hat I,\ j\in I.
	 $$
	 
	 Highest weight $U_q\,\widehat{\mathfrak{gl}}_{m|n}$-modules are admissible.


\begin{thebibliography}{0000000}
		
		\bibitem[AFS]{AFS}
		H. Awata, B. Feigin, J. Shiraishi,
		{\it Quantum algebraic approach to refined topological vertex},
		J. High Energy Phys. {\bf 2012} (2012), no. 3, 41--68
		
		\bibitem[BM]{BM} L. Bezerra and E. Mukhin,
		{\it Braid actions on quantum toroidal superalgebras}, arXiv:1912.08729
		
		\bibitem[BS]{BS}
		I. Burban and O. Schiffmann, 
		{\it On the Hall algebra of an elliptic curve I},
		Duke Math. J. {\bf 161} (2012), no. 7, 1171--1231
		
		\bibitem[DI]{DI}
		J. Ding and K. Iohara, 
		{\it Generalization of Drinfeld Quantum Affine Algebras},
			{Lett. Math. Phys.} {\bf 41} (1997), no. 2,  181--193
		
		\bibitem[FJM1]{FJM1} B. Feigin, M. Jimbo and E. Mukhin, 
		{\it Integrals of motion from quantum toroidal algebras}, 
		{\em J. Phys. A: Math. Theor.} {\bf 50} (2017) 464001
		
		\bibitem[FJM2]{FJM2} B. Feigin, M. Jimbo and E. Mukhin,
		{\it An evaluation homomorphism for quantum toroidal $\mathfrak{gl}(n)$ algebras}, arXiv:1709.01592v2 
		
		\bibitem[FJM3]{FJM3} B. Feigin, M. Jimbo and E. Mukhin, 
		{\it Towards trigonometric deformation of $\widehat{\mathfrak{sl}}_2$ coset VOA}, 
		arXiv:1811.02056v1
		
		\bibitem[FJMM1]{FJMM1} B. Feigin, M. Jimbo, T. Miwa and E. Mukhin, 
		{\it Representations of quantum toroidal $\mathfrak{gl}_N$}, 
		J. Algebra \textbf{380} (2013), 78--108
		
		\bibitem[FJMM2]{FJMM2} B. Feigin, M. Jimbo, T. Miwa and E. Mukhin, 
		{\it Branching rules for quantum toroidal $\mathfrak{gl}_N$}, 
		Adv. Math. \textbf{300} (2016), 229--274
		
		\bibitem[FJMM3]{FJMM3} B. Feigin, M. Jimbo, T. Miwa and E. Mukhin, 
		{\it Finite type modules and  Bethe ansatz for the quantum 
			toroidal $\mathfrak{gl}_1$}, Ann. Henri Poincar\'e \textbf{18} (2017), no. 8, 2543--2579
		
		\bibitem[FT]{FT}
		B. Feigin and A. Tsymbaliuk,
		{\it Heisenberg action in the equivariant
		K-theory of Hilbert schemes via Shuffle Algebra},
		Kyoto J. Math.  {\bf 51}  (2011),  no. 4, 831--854
		
		\bibitem[GKV]{GKV} V. Ginzburg, M. Kapranov, and E. Vasserot,
		{\it Langlands reciprocity for algebraic surfaces}, Math. Res. Lett. {\bf 2} (1995), no. 2, 147--160
		
		\bibitem[K1]{K1} T. Kojima,
		{\it A bosonization of $U_q(\widehat{\mathfrak{sl}}_{m|n})$}, 
		Comm. Math. Phys. {\bf 355} (2017),  no. 2, 603–-644
		
		\bibitem[K2]{K2} T. Kojima,
		{\it Commutation relations of vertex operators for $U_q(\widehat{\mathfrak{sl}}_{m|n})$}, 
		J. Math. Phys. {\bf 59} (2018),  no. 10, 101701, 37 pp. 
		
		\bibitem[KSU]{KSU} K. Kimura, J. Shiraishi, and J. Uchiyama,
		{\it A level-one representation of the quantum affine superalgebra  $U_q(\widehat{\mathfrak{sl}}(M+1|N+1))$}, 
		Comm. Math. Phys. {\bf 188} (1997),  no. 2, 367–-378
		
		\bibitem[KW]{KW} V. G. Kac and M. Wakimoto,
		{\it Integrable Highest Weight Modules over Affine Superalgebras and Appell's Function}, 
		Comm. Math. Phys. {\bf 215} (2001),  no. 3, 631–-682
		
		\bibitem[M1]{M1} K. Miki, {\it Toroidal braid group action and an automorphism
			of toroidal algebra $U_q\bigl(\mathfrak{sl}_{n+1,tor}\bigr)$ ($n\ge2$)},
		{Lett. Math. Phys.} {\bf 47} (1999), no. 4,  365--378
		
		\bibitem[M2]{M2} K. Miki,  {\it Toroidal and level $0$ $U_q'\widehat{sl_{n+1}}$ actions on $U_q\widehat{gl_{n+1}}$
			modules}, J. Math. Phys. {\bf 40} (1999), no. 6, 3191--3210
		
		
		\bibitem[N]{N} A. Negut,
		{\it The Shuffle Algebra Revisited},
		Int. Math. Res. Not. {\bf 2014} (2014), no. 22, 6242-–6275
		
		\bibitem[S]{S} Y. Saito,
		{\it Quantum toroidal algebras and their vertex representations}, 
		Publ. Res. Inst. Math. Sci. {\bf 34} (1998),  no. 2, 155–-177
		
		\bibitem[SV1]{SV1}
		O. Schiffmann and E. Vasserot,
		{\it The elliptic Hall algebra, Cherednik Hecke algebras and Macdonald
			polynomials}, Compos. Math. {\bf 147} (2011), no. 1, 188--234
		
		\bibitem[SV2]{SV2}
		O.~Schiffmann and E.~Vasserot,
		{\it The elliptic Hall algebra and the equivariant K-theory of the Hilbert scheme of $\mathbb{A}^2$},
		Duke Math. J. {\bf 162}  (2013),  no. 2, 279--366
		
		\bibitem[T1]{T1} A. Tsymbaliuk,
		{\it Quantum affine Gelfand-–Tsetlin bases and quantum toroidal algebra via K-theory of affine Laumon spaces}, 
		Sel. Math. New Ser. {\bf 16} (2010), no. 2,  173--200
		
		\bibitem[T2]{T2} A. Tsymbaliuk,
		{\it PBWD bases and shuffle algebra realizations for 
			$U_{\boldsymbol{v}}(L\mathfrak{sl}_n)$, $U_{\boldsymbol{v_1},\boldsymbol{v_2}}(L\mathfrak{sl}_n)$, $U_{\boldsymbol{v}}(L\mathfrak{sl}_{m|n})$}, arXiv:1808.09536
		
		\bibitem[VV1]{VV} M. Varagnolo and E. Vasserot, {\it Schur duality in the toroidal setting}, Comm. Math. Phys. {\bf 182} (1996), no. 2, 469--483
		
		\bibitem[Y]{Y} H. Yamane,
		{\it On defining relations of affine Lie superalgebras and affine quantized universal enveloping superalgebras}, 
		Publ. RIMS, Kyoto Univ. {\bf 35} (1999), 321-–390
		
			\bibitem[Z]{Z} Y. Zhang,
		{\it Comments on the Drinfeld realization of the quantum affine superalgebra $U_q[gl(m|n)^{(1)}]$ and its Hopf algebra structure}, 
		J. Phys. A: Math. Gen. {\bf 30} (1997), 8325-8335

		
	\end{thebibliography}
\end{document}